\documentclass[12pt]{amsart}

\usepackage{enumerate}

\input{diagrams.tex}
\diagramstyle[scriptlabels,height=8mm,width=8mm]

\def\pdfsyncstart{}
\def\pdfsyncstop{}
\def\bdi{\pdfsyncstop\begin{diagram}}
\def\edi{\end{diagram}\pdfsyncstart}

\theoremstyle{plain}

\newtheorem{thm}{Theorem}[section]
\newtheorem{cor}[thm]{Corollary}
\newtheorem{lem}[thm]{Lemma}
\newtheorem{prop}[thm]{Proposition}

\theoremstyle{definition}
\newtheorem{defi}[thm]{Definition}
\newtheorem{defis}[thm]{Definitions}
\newtheorem{conj}[thm]{Problem}
\newtheorem{conv}[thm]{Convention}
\newtheorem{nota}[thm]{Notation}
\newtheorem{rem}[thm]{Remark}
\newtheorem{rems}[thm]{Remarks}
\newtheorem{exa}[thm]{Example}
\newtheorem{exas}[thm]{Examples}
\newtheorem{sit}[thm]{}

\newcommand{\brem}{\begin{rem}}
\newcommand{\brems}{\begin{rems}}
\newcommand{\erem}{\end{rem}}
\newcommand{\erems}{\end{rems}}
\newcommand{\bexa}{\begin{exa}}
\newcommand{\bexas}{\begin{exas}}
\newcommand{\eexa}{\end{exa}}
\newcommand{\eexas}{\end{exas}}
\newcommand{\bdefi}{\begin{defi}}
\newcommand{\edefi}{\end{defi}}
\newcommand{\bdefis}{\begin{defis}}
\newcommand{\edefis}{\end{defis}}
\newcommand{\bcor}{\begin{cor}}
\newcommand{\ecor}{\end{cor}}
\newcommand{\blem}{\begin{lem}}
\newcommand{\elem}{\end{lem}}
\newcommand{\bconv}{\begin{conv}}
\newcommand{\econv}{\end{conv}}
\newcommand{\bconj}{\begin{conj}}
\newcommand{\econj}{\end{conj}}
\newcommand{\bprop}{\begin{prop}}
\newcommand{\eprop}{\end{prop}}
\newcommand{\bthm}{\begin{thm}}
\newcommand{\ethm}{\end{thm}}
\newcommand{\bnota}{\begin{nota}}
\newcommand{\enota}{\end{nota}}
\newcommand{\bsit}{\begin{sit}}
\newcommand{\esit}{\end{sit}}
\newcommand{\be}{\begin{equation}}
\newcommand{\ee}{\end{equation}}
\newcommand{\bproof}{\begin{proof}}
\newcommand{\eproof}{\end{proof}}
\def\ba{\begin{array}}
\def\ea{\end{array}}
\def\bea{\begin{eqnarray}}
\def\eea{\end{eqnarray}}

\def\bnum{\begin{enumerate}}
\def\enum{\end{enumerate}}
\newcommand{\no}{\noindent}

\def\lto{\longrightarrow}
\def\hto{\hookrightarrow}

\def\ext{{\rm ext}}
\def\sto{\rightsquigarrow}

\newcommand{\Spec}{\operatorname{Spec}}
\newcommand{\Div}{\operatorname{Div}}
\newcommand{\Pic}{\operatorname{Pic}}

\newcommand{\Frac}{\operatorname{Frac}}

\newcommand{\id}{\operatorname{id}}
\newcommand{\Der}{\operatorname{Der}}

\newcommand{\NS}{{\operatorname{NS}}}
\newcommand{\Aut}{{\operatorname{Aut}}}

\newcommand{\cy}{{\operatorname{Cycl}}}
\newcommand{\reg}{{\operatorname{reg}}}

\newcommand{\Ext}{{\operatorname{Ext}}}

\def\supp{{\rm supp\,}}

\renewcommand{\div}{{\operatorname{div}}}

\def\fF{{\mathfrak F}}

\def\fm{{\mathfrak m}}

\def\cB{{\mathcal B}}
\def\cC{{\mathcal C}}
\def\cD{{\mathcal D}}
\def\cE{{\mathcal E}}
\def\cF{{\mathcal F}}

\def\cL{{\mathcal L}}

\def\cN{{\mathcal N}}
\def\cO{{\mathcal O}}
\def\cR{{\mathcal R}}

\def\cX{{\mathcal X}}
\def\cV{{\mathcal V}}

\def\cZ{{\mathcal Z}}

\def\bB{{\bar B}}
\def\bC{{\bar C}}
\def\bD{{\bar D}}

\def\bV{{\bar V}}
\def\bX{{\bar X}}

\def\bF{{\bar F}}

\def\tC{{\tilde C}}
\def\tO{{\tilde O}}
\def\tD{{\tilde D}}
\def\tR{{\tilde R}}
\def\tB{{\tilde B}}
\def\tV{{\tilde V}}
\def\tX{{\tilde X}}

\def\tgamma{{\tilde\gamma}}

\newcommand{\la}{\label}

\newcommand{\A}{{\mathbb A}}
\newcommand{\Di}{{\mathbb D}}
\newcommand{\PP}{{\mathbb P}}

\newcommand{\C}{{\mathbb C}}
\newcommand{\Q}{{\mathbb Q}}
\newcommand{\Z}{{\mathbb Z}}
\newcommand{\N}{{\mathbb N}}
\newcommand{\T}{{\mathbb T}}

\newcommand{\G}{{\Gamma}}

\newcommand{\p}{{\partial}}

\def\bO{{\bar O}}
\def\quot{/\hskip-3pt/}

\def\reg{{reg}}
\def\hF{{\hat F}}

\newcommand{\nlin}{\unitlength1mm\begin{picture}(0,9.25)
                     \put(0,0.75){\line(0,1){8.5}}
                    \end{picture}}

\newcommand{\vlin}[1]{\hspace{0.75mm}\unitlength1mm\begin{picture}
(#1,0)
                     \put(0,0){\line(1,0){#1}}
                    \end{picture}\hspace{0.75mm}\rule[-3mm]{0mm}{4mm}}

\def\llin{\vlin{11.5}}
\newcommand{\lin}{\vlin{8.5}}

\newcommand{\co}[1]{\unitlength1mm\begin{picture}(0,8)
  \put(0,0){\circle{1.5}}
  \put(0,3){\makebox(0,5)[b]{$#1$}}
                    \end{picture}}
\newcommand{\mybox}{\unitlength1mm\begin{picture}(0,1.5)
  \put(-0.75,-0.75){\line(0,1){1.5}}
  \put(-0.75,-0.75){\line(1,0){1.5}}
  \put(0.75,0.75){\line(0,-1){1.5}}
  \put(0.75,0.75){\line(-1,0){1.5}}
  \end{picture}}
\newcommand{\boxo}[1]{\unitlength1mm\begin{picture}(0,8)
  \put(0,0){\mybox}
  \put(0,3){\makebox(0,5)[b]{$#1$}}
                    \end{picture}}
\newcommand{\xbox}{\unitlength1mm\begin{picture}(0,1.5)
  \put(0,0){$\mybox$}
  \put(-0.75,0){\line(1,0){1.5}}
  \put(0,-0.75){\line(0,1){1.5}}
  \end{picture}}
\newcommand{\xboxo}[1]{\unitlength1mm\begin{picture}(0,8)
  \put(0,0){\xbox}
  \put(0,3){\makebox(0,5)[b]{$#1$}}
                    \end{picture}}
\newcommand{\cou}[2]{\unitlength1mm\begin{picture}(0,8)
  \put(0,0){\circle{1.5}}
  \put(0,3){\makebox(0,5)[b]{$#1$}}
  \put(0,-7){\makebox(0,4)[t]{$#2$}}
    \end{picture}
    \rule[-7mm]{0mm}{7mm}}

\newcommand{\boxou}[2]{\unitlength1mm\begin{picture}(0,8)
  \put(0,0){\mybox}
  \put(0,3){\makebox(0,5)[b]{$#1$}}
  \put(0,-7){\makebox(0,4)[t]{$#2$}}
    \end{picture}
    \rule[-7mm]{0mm}{7mm}}

\newcommand{\crl}[2]{\unitlength1mm\begin{picture}(0,8)
  \put(0,0){\circle{1.5}}
  \put(-5,0){\makebox(0,5)[b]{$#1$}}
 \put(5,0){\makebox(0,5)[b]{$#2$}}
    \end{picture}
    \rule[-7mm]{0mm}{7mm}}

\newcommand{\cshiftup}[2]{\unitlength1mm\begin{picture}(0,9.25)
                     \put(0,10){\crl{#1}{#2}}
                    \end{picture}}

\newcommand{\xbrl}[2]{\unitlength1mm\begin{picture}(0,8)
  \put(0,0){\xbox}
  \put(-5,0){\makebox(0,5)[b]{$#1$}}
 \put(5,0){\makebox(0,5)[b]{$#2$}}
    \end{picture}
    \rule[-7mm]{0mm}{7mm}}

\newcommand{\boxrl}[2]{\unitlength1mm\begin{picture}(0,8)
  \put(0,0){\mybox}
  \put(-5,0){\makebox(0,5)[b]{$#1$}}
 \put(5,0){\makebox(0,5)[b]{$#2$}}
    \end{picture}
    \rule[-7mm]{0mm}{7mm}}

\newcommand{\xbshiftup}[2]{\unitlength1mm\begin{picture}(0,9.25)
                     \put(0,10){\xbrl{#1}{#2}}
                    \end{picture}}

\newcommand{\boxshiftup}[2]{\unitlength1mm\begin{picture}(0,9.25)
                     \put(0,10){\boxrl{#1}{#2}}
                    \end{picture}}

\newcommand{\rbh}[1]{\raisebox{1mm}[1mm][-1mm]{#1}}

\title{Uniqueness of $\C^*$- and $\C_+$-actions on Gizatullin
surfaces}

\author{Hubert Flenner}
\address{Fakult\"at f\"ur Mathematik,
Ruhr Universit\"at Bochum, Geb.\ NA 2/72, Universit\"ats\-str.\
150, 44780 Bochum, Germany}
\email{Hubert.Flenner@ruhr-uni-bochum.de}

\author{Shulim Kaliman}
\address{Department of Mathematics,
University of Miami, Coral Gables, FL  33124, U.S.A.}
\email{kaliman@math.miami.edu}

\author{Mikhail Zaidenberg}
\address{Universit\'e
Grenoble I, Institut Fourier, UMR 5582 CNRS-UJF, BP 74, 38402
St.\ Martin d'H\`eres c\'edex, France}
\email{zaidenbe@ujf-grenoble.fr}

\thanks{
{\bf Acknowledgements:} This research was done during a visit of
the first and the second authors at the Institut Fourier, Grenoble,
of the third one at the Ruhr University at Bochum,
 and of all three authors at
the Max-Planck-Institute of Mathematics, Bonn.
They thank these institutions for the
generous support and excellent working conditions.}

\thanks{
\mbox{\hspace{11pt}}{\it 1991 Mathematics Subject
Classification}:
14R05, 14R20, 14J50.\\
\mbox{\hspace{11pt}}{\it Key words}: $\C^*$-action, $\C_+$-action,
affine surface}

\date{}

\begin{document}

\begin{abstract}
A {\it Gizatullin surface} is a normal affine surface $V$ over $\C$,
which can be completed by a zigzag; that is, by a linear chain of
smooth rational curves. In this paper we deal with the question
of uniqueness of $\C^*$-actions and $\A^1$-fibrations on such a
surface  $V$
up to automorphisms. The latter
fibrations are in one to one correspondence with $\C_+$-actions on
$V$ considered up to a ``speed change".

Non-Gizatullin surfaces are known to admit at most one
$\A^1$-fibration $V\to S$ up to an isomorphism of the base $S$.
Moreover an effective $\C^{*}$-action on them, if it does exist,
is unique up to conjugation and inversion $t\mapsto t^{-1}$ of
$\C^*$. Obviously uniqueness of $\C^*$-actions fails for affine
toric surfaces; however we show in this case that there are at
most two conjugacy classes of $\A^1$-fibrations. There is a
further interesting family of non-toric Gizatullin surfaces,
called the Danilov-Gizatullin surfaces, where there are in general
several conjugacy classes of $\C^*$-actions and $\A^1$-fibrations,
see e.g., \cite{FKZ1}.

In the present paper we obtain a criterion
as to when $\A^1$-fibrations of
Gizatullin surfaces are conjugate up to an automorphism of $V$ and the base $S$.
We exhibit as well a large subclasses of Gizatullin $\C^{*}$-surfaces
for which a $\C^*$-action
is essentially unique and
for which
there are at most two conjugacy classes
of $\A^1$-fibrations over $\A^1$.
\end{abstract}

\maketitle

\tableofcontents

\section*{Introduction}

Let $V$  be a normal affine surface admitting an effective action
of the group  $\C^*$. It is a natural question as to when any two
such actions on $V$ are conjugate in the automorphism group $\Aut
(V)$. Similarly, given an  $\C_+$-action on $V$ one may ask
whether its associated $\A^1$-fibration $V\to S$ is unique up to
conjugation; that is,  up to an automorphism of $V$ and an
isomorphism of the the base $S$.

Recall \cite{FKZ2} that a {\it Gizatullin surface} is a normal
affine surface completable by a {\it zigzag} that is, by a linear
chain of smooth rational curves. The uniqueness of $\C^*$-actions
on normal affine surfaces, up to conjugation and inversion, is
known to hold for all non-Gizatullin surfaces (see \cite{Be} for
the smooth case, \cite[Theorem 3.3]{FlZa3} for the general one).
Similarly in these cases there is at most one $\A^1$-fibration
$V\to S$ over an affine base up to an isomorphism of $S$, so any two
$\C_+$-actions define the same $\A^1$-fibration. However
uniqueness fails for every affine toric surface, which admits a
sequence of pairwise non-conjugate $\C^*$-actions.

Another important class of counterexamples is provided by the {\em
Danilov-Gizatullin surfaces}. By definition such a surface is the
complement of an ample section say $S$ in a Hirzebruch surface
$\Sigma_n$. A surprising theorem established in
\cite{DaGi}\footnote{See \cite[Corollary 4.8]{CNR} for an
alternative approach.} says that the isomorphism type of such a
surface $V_{k+1}=\Sigma_n\setminus S$ depends only on $k:=S^2-1$ and
neither on $n$ nor on $S$.
Answering our question,
Peter Russell observed that
the Danilov-Gizatullin theorem  actually
provides $k$
pairwise non-conjugate
$\C^*$-actions on $V_{k+1}$. We reproved  in
\cite[5.3]{FKZ2} this result showing moreover
that these $k$ $\C^*$-actions exhaust all
$\C^*$-actions on $V_{k+1}$ up to conjugation.
At least half of them stay
non-conjugate up to inversion
in $\C^*$. Moreover by \cite[5.16]{FKZ2} in this case
there are at least $\lfloor \frac{k+1}{2}\rfloor$
different conjugacy classes of $\A^1$-fibrations.

Let us recall that every Gizatullin surface $V\not\cong
\A^1\times\C^*$ can be completed by a {\em standard} zigzag \be
\label{izigzag} \cou{C_0}{0}\lin\cou{C_1}{0}\lin\cou{C_2}{w_2}
\lin\ldots\lin\cou{C_n}{w_n}\quad, \ee with $w_i=C_i^2\le
-2\;\forall i\ge 2$. Although this completion is not unique the
sequence of weights $(w_2,\ldots, w_n)$ is up to reversion an
invariant of $V$ \cite{Gi}, cf. also \cite{Du1,FKZ2}.

The linear system $|C_0|$ provides a $\PP^1$-fibration
$\Phi_0:\bV\to \PP^1$, which restricts to an $\A^1$-fibration
$\Phi_0:V\to \A^1$ (similarly, reversing the zigzag gives a second
$\A^1$-fibration  $\Phi_0^\vee:V\to \A^1$).
This $\PP^1$-fibration lifts to the minimal resolution of
singularities $\tV$ of $\bar V$. Our results are formulated in
terms of the so called {\em extended boundary divisor}
$$
D_\ext:=C_0+C_1+\tilde\Phi_0^{-1}(0)\subseteq \tV
$$
considered in \cite{Gi,
Du1, FKZ2}, where $\tilde\Phi_0$ is the induced fibration. Its
structure is well known, see Proposition \ref{ext}. We introduce
{\em rigid}\/ and {\em distinguished}\/ extended
divisors that are characterized by their weighted dual graph, see
\ref{dist} and \ref{rig} for details. The main result of the paper
(see Theorem \ref{main}) can be stated as follows.

\bthm\label{01} Let $V$ be a Gizatullin surface whose extended
divisor $D_\ext$ is distinguished and rigid. Then $\Phi_0$ and
$\Phi_0^\vee$ are up to conjugation the only $\A^1$-fibrations $V\to
\A^1$. \ethm

In the special case of surfaces $xy=p(z)$ in $\A^3$, this result was
obtained in terms of locally nilpotent derivations by Daigle
\cite{Dai} and Makar-Limanov \cite{ML2}.

Our approach has important applications to the classification of
$\C^*$-actions on $V$. In \cite{FKZ2} we conjectured that among
smooth affine $\C^*$-surfaces, the toric surfaces and the
Danilov-Gizatullin surfaces are the only exceptions to uniqueness of
a $\C^*$-action. In Theorem \ref{MT} below we confirm this
conjecture in the particular case of Gizatullin surfaces with a
rigid extended divisor. Recall \cite{FlZa1} that every normal affine
surface $V$ with a hyperbolic $\C^*$-action admits a {\it DPD
presentation} $V=\Spec A_0[D_+,D_-]$, where $D_+,D_-$ are two
$\Q$-divisors on the smooth affine curve $C=\Spec\,A_0$ with
$D_++D_-\le 0$, and $A_0$ is the ring of invariants; see
Section 3.1 for details. For a Gizatullin $\C^*$-surface $V$ one has
\cite{FlZa2}: $A_0=\C[t]$, and each of the fractional parts
$\{D_\pm\}=D_\pm-\lfloor D_\pm\rfloor$ is concentrated on at most
one point $\{p_\pm\}$. To formulate our second main result we
consider the following 3 conditions on $D_+$, $D_-$.

\begin{enumerate}
\item[($\alpha_+$)] $\supp \{D_+\}\cup \supp \{D_-\}$ is empty or
consists of one point, say, $p$ satisfying either
$D_+(p)+D_-(p)=0\,$  or
$$
D_+(p)+D_-(p)\le -\max\left(\frac{1}{{m^+}^2},\,
\frac{1}{{m^-}^2}\right)\,,$$
where $\pm m^\pm$ is the minimal
positive integer such that $m^\pm D_\pm(p)\in\Z$.

\item[($\alpha_*$)] $\supp \{D_+\}\cup \supp \{D_-\}$ is empty or
consists of one point  $p$, where $D_+(p) + D_-(p) \le -1$
or both fractional parts $\{D_+(p)\}$, $\{ D_-(p)\}$ are nonzero.

\item[($\beta$)] $\supp \{D_+\}=\{p_+\}$ and $\supp
\{D_-\}=\{p_-\}$ for two different points $p_+,p_-$, where
$D_+(p_+) + D_-(p_+)\le-1$ and $D_+(p_-) + D_-(p_-)\le-1.$
\enum

\bthm\label{MT}
For a non-toric normal Gizatullin
$\C^*$-surface $V=\Spec \C[t][D_+,D_-]$ the following hold.

1. If $(\alpha_*$) or $(\beta$) is fulfilled then the
$\C^*$-action on $V$ is unique up to conjugation in the
automorphism group of $V$ and up to inversion $\lambda\mapsto
\lambda^{-1}$ in $\C^*$. Moreover the given $\C^*$-action is
conjugate to its inverse if and only if for a suitable
automorphism $\psi\in\Aut (\A^1)$ \be\label{MTiff} \psi^*(D_+)-D_-
\mbox{ is integral and}
 \quad
\psi^*(D_++D_-)=D_++D_-\,.
\ee

2.  If $(\alpha_+$) or $(\beta)$ holds then up to conjugation
there are at most two conjugacy classes of $\A^1$-fibrations
$V\to\A^1$. There is only one such conjugacy class if and only if
(\ref{MTiff}) is fulfilled for some $\psi\in\Aut (\A^1)$. \ethm

We notice that for smooth non-toric Gizatullin $\C^*$-surfaces
this proves uniqueness of $\C^*$-actions up to conjugation and
inversion unless the weights $w_i$ in the boundary zigzag
(\ref{izigzag}) satisfy $w_i=-2$ $\forall i\ne s$ for some $s$ in
the range $2\le s\le n$. In a forthcoming paper we will show that in the
latter case there is a deformation family of pairwise non-conjugate
$\C^*$-actions on $V$. Consequently, for smooth Gizatullin
$\C^*$-surfaces the sufficient conditions in Theorem \ref{MT} are
also necessary.

Let us survey the content of the different sections.
In Section 1.1 we review some standard facts on Gizatullin surfaces
and describe in Section 1.2 their extended divisors. After some
preparations in 1.3 we treat in Section 1.4 families of completions
of a given Gizatullin surface by zigzags. The main result here is
the triviality criterion \ref{propdist}, which provides one of the
basic tools in the proof of Theorem \ref{01}. In Section 2 the
possible degenerations of extended divisors in such families are
studied. The main result here is Theorem \ref{maincrit}, which gives
a criterion for when the extended divisor is rigid, i.e.\ stays
constant in a family.

In Section 3 we translate these conditions into the language of
DPD presentations. First we recall  the description of standard
equivariant completions of Gizatullin $\C^*$-surfaces in terms of
a DPD presentation according to \cite{FKZ2}. In Theorem
\ref{nonspec} we give the required criterion
for the extended divisor $D_\ext$ to be distinguished and rigid.

One of our main technical tools is the so called {\em
reconstruction space}. Roughly speaking, the latter forms  a
moduli space for the completions of a given normal surface. In
Section 4 we show that this moduli space exists and is isomorphic
to an affine space, see Corollary \ref{rec.8}. This fact is a
basic ingredient in the proofs of Theorems \ref{01} and \ref{MT}
in the final Section 5.

\section{Gizatullin surfaces}

\subsection{Standard completions of Gizatullin surfaces}
Let us recall the notion of a standard zigzag
\cite{FKZ1}.

\bsit Let $X$ be a complete normal algebraic surface. By a {\it
zigzag} on $X$ we mean an SNC divisor\footnote{I.e. a simple normal crossing divisor.}
$D$ with rational components
contained in the smooth part $X_\reg$, which has a linear dual graph
\be\label{zigzag} \G_D:\quad\quad
\cou{C_0}{w_0}\lin\cou{C_1}{w_1}\lin\ldots\lin\cou{C_n}{w_n}\quad,
\ee where $w_0,\ldots, w_n$ are the weights of $\G_D$. We abbreviate
this chain by $[[w_0,\ldots, w_n]]$. We also write $[[\ldots,
(w)_k,\ldots ]]$ if a weight $w$ occurs at $k$ consecutive places.
Note that the intersection matrix of a zigzag has at most one
positive eigenvalue by the Hodge index theorem. We recall the
following notion. \esit

\bdefi\label{sst} (\cite[Definition 2.13 and Lemma 2.17]{FKZ1}) A
zigzag $D$ is called {\em standard}  if its dual graph $\G_D$ is
one of
\be \label{standardzigzag} [[0]]\,,\quad [[0,0]]\,, \quad
[[0,0,0]] \,\quad\mbox{or}\quad
  [[0,0,w_2,\ldots,w_n]],\quad \mbox{where}
\quad n\ge 2,\, w_j\le-2\,\,\,\forall j. \ee
A linear chain $\Gamma$ is said to be {\em semistandard}
if it is either standard or one of
\be \label{sstandard}
 [[0,w_1,w_2,\ldots,w_n]],\quad [[0,w_1,0]] \quad \mbox{where}
\quad m\in \Z,\;  n\ge 1,\, w_j\le-2\,\,\,\forall j. \ee We note
that a standard zigzag $[[0,0,w_2,\ldots,w_n]]$ is unique in its
birational class up to reversion \be\la{reverse}
[[0,0,w_2,\ldots,w_n]]\rightsquigarrow [[0,0,w_n,\ldots,w_2]]\,, \ee
see Corollary 3.33 in \cite{FKZ1}.  A zigzag is called {\em
symmetric} if it coincides with its reversed zigzag.   \edefi

By definition a {\it Gizatullin surface} is a normal affine surface
$V$  which admits a completion $(\bV, D)$ with a zigzag $D$. Such a
completion is called {\em (semi)standard} if $D$ has this property.
We need  the following facts.

\blem
\label{equivariant} For a Gizatullin surface $V$ the
following hold.
\bnum[(a)]
\item (\cite[Corollary 3.36]{DaGi, Du1,
FKZ1}) $V$ admits a standard completion $(\bV,D)$.
\item (\cite[Theorem 2.9(b)]{FKZ2}) If a torus $\T=(\C^*)^m$
acts on $V$ then $V$ admits an
equivariant standard completion, which is unique up to reversing
the boundary zigzag.
\item
(\cite[Theorem 2.9(a) and Remark 2.10(1)]{FKZ2})
If $\C_+$ acts on $V$ then $V$ admits an
equivariant {\it semistandard} completion.
\enum
\elem

\bsit\label{reversion} The reversion of a zigzag, regarded as a
birational transformation of the weighted dual graph, admits the
following factorization \cite{FKZ1}. Given
$[[0,0,w_2,\ldots,w_n]]$ we can successively move the pair of zeros
to the right
$$
[[0,0,w_2,\ldots,w_n]]\sto [[w_2, 0,0,w_3,\ldots,w_n]]
\sto \ldots \sto [[w_2,\ldots,w_n, 0,0]]
$$
by a sequence of {\em inner elementary transformations}\footnote{
By an inner elementary transformation of a weighted graph we mean
blowing up at an edge incident to a $0$-vertex of degree $2$ and
blowing down the image of this vertex.}, see Example 2.11(2) in
\cite{FKZ1}. The corresponding birational transformation
$[[0,0,w_2,\ldots,w_n]]\sto [[w_2,\ldots,w_n, 0,0]]$ is
non-trivial unless our standard graph is one of $[[0]]$, $[[0,0]]$
or $[[0,0,0]]$.

If $(\bV,D)$ is a standard completion of a Gizatullin surface $V$,
then reversing the zigzag $D$ by a sequence of inner elementary
transformations as explained above we obtain from $(\bV, D)$ a new
completion  $(\bV^\vee, D^\vee)$, which we call the {\em reverse
standard completion}. It is uniquely determined by $(\bV, D)$.
Note that even in the case where the zigzag $D$ is symmetric
with dual graph $\ne [[0]]$, $[[0,0]]$, $[[0,0,0]]$, this reverse
completion $(\bV^\vee, D^\vee)$ is not isomorphic to $(\bV, D)$
under an isomorphism fixing pointwise the affine part $V$. \esit

\subsection{Extended divisors of Gizatullin surfaces}
\bsit\la{unisit} Let $V$ be a
Gizatullin surface and
$(\bV, D)$ be a completion of $V$ by a standard zigzag
$[[0,0,w_2,\ldots, w_n]]$ with $n\ge 2$  and $w_i\le -2$ $\forall
i$. We write
$$
D=C_0+\ldots+C_n\,,
$$
where the irreducible components $C_i$ are enumerated as in
(\ref{zigzag}). We consider the minimal resolutions of
singularities $V'$, $(\tV,D)$ of $V$  and $(\bV, D)$,
respectively.

Since $C_0^2=C_1^2=0$, the linear systems $|C_0|$ and $|C_1|$ define
a morphism $\Phi=\Phi_0\times\Phi_1:\tV\to \PP^1\times\PP^1$ with
$\Phi_i=\Phi_{|C_i|}$, $i=0,1$. We call it the {\it standard
morphism} associated to the standard completion $(\bV, D)$ of $V$.
Similarly $\Phi_0$ is referred to as the {\it standard
$\PP^1$-fibration} of $(\bV, D)$.

We note that $C_1$ is a section of $\Phi_0$ and so the
restriction $\Phi_0|V':V'\to\PP^1$ is an $\A^1$-fibration. We can
choose the coordinates on $\PP^1=\C\cup \{\infty\}$ in such a way
that
$$C_0=\Phi_0^{-1}(\infty)\,, \qquad \Phi(C_1)= \PP^1\times
\{\infty\} \quad\mbox{and}\quad C_2\cup\ldots \cup C_n\subseteq
\Phi_0^{-1}(0) \,.$$ The standard morphism $\Phi$ contracts the
curves $C_i$ for $i\ge 3$ and does not contract $C_0,C_1,C_2$. By
abuse of notation we denote the images of $C_0,C_1,C_2$ in
$\PP^1\times\PP^1$ by the same letters. The divisor
$D_\ext:=C_0\cup C_1\cup \Phi_0^{-1}(0)$ is called the {\it
extended divisor}.
\esit

\brem\label{unirem}
1. The  dual graph of $D_\ext$ is linear if and only if $V$ is
toric \cite[Lemma 2.20]{FKZ2}.

2. If $V$ carries a $\C^*$-action then we can find an equivariant
standard completion  $(\bV, D)$, see Lemma \ref{equivariant}(b).
Since the minimal resolution of singularities is also equivariant,
so are $(\tV,D)$ and $\Phi$ with respect to a suitable $\C^*$-action
on $\PP^1\times\PP^1$, and the divisor $D_\ext$ is invariant under
the  $\C^*$-action on $\tV$. For $\C^*$-surfaces this divisor was
studied systematically in \cite {FKZ2}.

3. The morphism $\Phi=\Phi_0\times\Phi_1$ contracts $C_3\cup\ldots
\cup C_n$, in particular it contracts all exceptional curves in the
resolution $V'\to V$, whence descends to a morphism
$\bar\Phi=\bar\Phi_0\times\bar\Phi_1:\bar V\to \PP^1\times\PP^1.$
We also call $\Phi$ the standard morphism of $(\bV, D)$ and
$\bar\Phi_0$
the standard $\PP^1$-fibration.
\erem

We recall the following fact, see \cite[Lemma 2.19]{FKZ2}.

\blem\label{unilem}
With the notation as in \ref{unisit},
$\Phi$ is birational and induces an isomorphism
$\tV\backslash \Phi_0^{-1}(0)\cong
(\PP^1\backslash\{0\})\times \PP^1$.
In particular, $D_{(e)}:= \Phi_0^{-1}(0)$ is the only
possible degenerate fiber of the
$\PP^1$-fibration $\Phi_0 : \tV\to \PP^1$.
\elem

To exhibit the structure of this extended divisor let us recall
some notation from \cite{FKZ2}.

\bsit\la{toricsit} For a primitive $d$th root of unity $\zeta$ and
$0\le e< d$ with $\gcd(e,d)=1$ \footnote{In the case $d=1$ this
forces $(d,e)=(1,0)$.}
 the cyclic group $\Z_d=\langle \zeta
\rangle$ acts on $\A^2$ via $\zeta.(x,y)=(\zeta x, \zeta^e y)$. The
quotient $V_{d,e}=\A^2\quot \Z_d$ is a normal affine toric surface.
Moreover, any such surface different from\footnote{Hereafter
$\A^1_*=\A^1\setminus\{0\}$.}  $\A^1_*\times\A^1_*$ and
$\A^1_*\times\A^1$ arises in this way. Singularities analytically
isomorphic to the singular point of $V_{d,e}$ are called cyclic
quotient singularities of type $(d,e)$. \esit

\bsit\label{not1} We abbreviate by a box $\rbh{ \mybox\,}$ with
rational weight $e/m$, where $0< e<m$ and $\gcd(m,e)=1$, the
weighted linear graph \be\label{gr1} \cou{C_1}{-k_1}
\lin\ldots\lin\cou{C_n}{-k_n} \qquad=\qquad\boxo{e/m}\ee with
$k_1,\ldots,k_n\ge 2$, where
$$m/ e=[k_1,\ldots, k_n]=k_1-\frac{1}{k_2-\frac{1}{\ddots
-\frac{1}{k_n}}}\quad.$$ A chain of rational curves $(C_i)$ on a
smooth surface with dual graph (\ref{gr1}) contracts to a cyclic
quotient singularity of type $(m,e)$ \cite{Hi}. It is convenient
to introduce the weighted box $\quad\boxo{0}\qquad$ for the empty
chain. Given extra curves $E,\,F$ we also abbreviate
\be\label{reversed chain} \co{E}{}\lin
\co{C_1}\lin\ldots\lin\co{C_n} \qquad
=\qquad\co{E}\llin\boxo{e/m}\qquad = \qquad
\boxo{(e/m)^*}\llin\co{E}\quad\,\, \ee and \be\label{reversed
chain1} \co{C_1} \lin\ldots\lin\co{C_n} \lin\co{F} \qquad=\qquad
\boxo{e/m}\llin\co{F}\qquad
=\qquad\co{F}\llin\boxo{(e/m)^*}\quad.\ee The orientation of the
chain of curves $(C_i)_i$ in (\ref{gr1}) plays an important role.
Indeed
 $[k_n,\ldots,k_1]=m/e'$, where $0< e'<m$,
$ee'\equiv 1 \pmod m$, and the box $\rbh{ \mybox}\,$ marked with
$(e/m)^*:=e'/m$ corresponds to the reversed chain in (\ref{gr1}),
see e.g., \cite{Rus}. The chain $[[(-2)_m]]$ will be abbreviated
by \;\;$\mybox\;\; A_m$. \esit

\bdefi\label{eboundary.4} A {\em feather} $\fF$ is a linear chain of
smooth rational curves with dual graph \be\label{picfeather}
\fF:\qquad\qquad \co{B}\llin\boxo{e/m} \quad,\ee where $B$ has
self-intersection $\le -1$ and $e,m$ are as before, cf.
(\ref{reversed chain}). Note that the box does not contain a
$(-1)$-curve; it can also be empty. The curve $B$ will be called the
{\em bridge curve}.

A collection of feathers
$\{\fF_\rho\}$ consists of feathers $\fF_\rho$, $\rho=1,\ldots
,r$, which are pairwise disjoint. Such a collection will be
denoted by a plus box $\rbh{ \xbox}\;$.
We say that a collection $\{\fF_\rho\}$ is attached
to a curve $C_i$ in a chain (\ref{zigzag}) if the bridge curves
$B_\rho$ meet $C_i$ in pairwise distinct points and all the
feathers $\fF_\rho$
are disjoint with the curves $C_j$ for $j\ne i$. In a
diagram we write in brief
$$
\co{C_i}\lin\xboxo{\{\fF_\rho\}} \qquad\mbox{or, in the case of a
single feather,}\qquad \co{C_i}\lin\boxo{\fF}\quad.
$$
We often draw this diagram vertically, with the same meaning.
\edefi

An  {\em $A_k$-feather} $\quad \co{B}\lin \boxo{A_k}\quad $
represents the contractible linear chain $[[-1,(-2)_{k}]]$. Thus the
$A_0$-feather represents a single $(-1)$-curve $B$, while the box is
empty.

Let us further exhibit the structure
of the extended divisor of a
Gizatullin
surface according to \cite{Du1}.

\bprop\label{ext} Let $(\tV,D)$ be a minimal SNC completion of the
minimal resolution of singularities of a Gizatullin surface $V$,
where $D=C_0+\ldots + C_n$ is a zigzag as in (\ref{zigzag}). Then
the extended divisor $D_\ext$ has dual graph
\bigskip

\be\label{extended} D_\ext: \qquad\quad\cou{0}{C_0}\llin
\cou{0}{C_1}\lin\cou{}{C_2} \nlin\xbshiftup{}{\,\,\{\fF_{2j}\}}
\llin \ldots\lin\cou{}{C_i}\nlin\xbshiftup{}{\,\,\{\fF_{ij}\}}
\llin\ldots\llin \cou{}{C_n} \nlin\xbshiftup{}{\,\,\{ \fF_{nj}\}}
\qquad\qquad, \ee where $\fF_{ij}$ ($1\le j\le r_i$) is a collection
of feathers attached to the curve $C_i$, $i\ge 2$. Moreover the
surface $\tV$ is obtained from $ \PP^1\times\PP^1$ by a sequence of
blowups with centers in the images of the components $C_i$, $i\ge
2$. \eprop

\bproof A proof can be found (using different notation) in
\cite{Du1}. For the convenience of the reader we provide a short
argument. First we note that $D_ {(e)}=\Phi_0^{-1}(0)\subseteq\tV$
is a tree of rational curves, since it is the blowup of a fiber
$C_2=\{0\}\times \PP^1\subseteq \PP^1\times\PP^1$. Let $ \fF_{ij}$,
$j=1,\ldots, r_i$,  be the connected components of $D_\ext\ominus
C_i$ that do not contain components of $D$. Every such connected
component contains a unique curve $B_{ij}$, which meets $C_i$. The
divisor $R_{ij}:=\fF_{ij}\ominus B_{ij}$ is then disjoint from $D$.
Since $V$ is affine, $R_{ij}$ contracts to a point in $V$. Hence it
is the exceptional divisor of a minimal resolution of a singular
point of $V$ and so its dual graph contains no linear $(-1)$-curve.
On the other hand, the divisor $D_{(e)}$ contracts
to $C_2$. We claim that the dual graph of
$\fF_{ij}$ contains no branch point, and its
end point $B_{ij}$ is the only possible $(-1)$-curve in $\fF_{ij}$.
Let us check this claim by induction on the number of blowdowns in
the contraction of $D_{(e)}$ to $C_2$, or rather of blowups when
growing $D_{(e)}$ starting from $C_2$. Indeed along this process

\begin{enumerate}[$\bullet$] \item
the image, say $D'$, of the chain $D=C_0+\ldots+C_n$ under an
intermediate blowdown is again a linear chain, \item the image
$\fF_{ij}'$ of $\fF_{ij}$ is either empty or a connected component
of $D_\ext'\ominus D'$, where $D_\ext'$ is the image of $D_\ext$.
Moreover \item if $\fF_{ij}'\neq\emptyset$ then it still contains
just one neighbor say $B_{ij}'$ of $D'$ in $D_\ext'$, \item
$R_{ij}':=\fF_{ij}'\ominus B_{ij}'$ is either empty or a minimal
resolution of a singular point with a linear dual graph, \item
$B_{ij}'$ is at most linear vertex in the dual graph of $D_\ext'$,
and the only possible $(-1)$-curve in $\fF_{ij}'$. \end{enumerate}

The next blowup must be done at a point of $D'$ (which is either a
smooth point of $D_\ext'$, or a double point of $D_\ext'$). Indeed
otherwise it would be done at a point of $\fF_{ij}'\ominus D'$,
and then clearly $R_{ij}$ cannot be minimal i.e., it would contain
a $(-1)$-curve, which is impossible. Thus all the properties
mentioned above are preserved under this blowup.

This implies that $\fF_{ij}$ is a linear feather of the form
$$\fF_{ij}:\quad\quad\co{B_{ij}}\lin\boxo{R_{ij}}\quad\,,$$
which yields the desired form of $D_\ext$, and also the last
assertion.
\eproof

\brem\label{smooth} The collection of linear chains $R_{ij}$
corresponds to the minimal resolution of singularities of $V$. So
$V$ has at most cyclic quotient singularities, cf. \cite[Ch. 3,
Lemma 1.4.4(1)]{Miy}. Moreover $V$ is smooth if and only if the
collection $R_{ij}$ is empty,  if and only if every feather
$\fF_{ij}$ reduces to a single bridge curve $B_{ij}$.\erem

\subsection{Simultaneous contractions}
The following lemma is a standard fact in surface theory.

\blem\label{rec.1}
For a smooth rational surface  $X$
and a smooth rational
curve $C$ on $X$ with $C^2=0$,
we have
$$
H^0(X,\cO_X(C))\cong \C^2\,,\qquad
H^i(X,\cO_X(C))=0
\mbox{ for } i\ge 1\,.
$$
Moreover the linear system $|C|$ is base point free and defines a
 $\PP^1$-fibration   $\Phi_{|C|}:X\to \PP^1$. \elem

A relative version of this result is as follows.

\begin{lem}\label{rec.2}
Let $f: \cX \to S$ be a smooth family of rational
surfaces over a quasiprojective scheme $S$
with $\Pic(S)=0$,
and let $\cC$ be an $S$-flat
divisor in $\cX$ such that the fibers
$\cC_s:=f^{-1}(s) \cap\cC$
are smooth rational curves
of self-intersection 0 in $\cX_s:=f^{-1}(s)$.
Suppose that
$R \subset \cX$ is a section of $f$
disjoint from $\cC$. Then there exists a morphism
$\varphi : \cX\to \PP^1$ such that $\varphi^{*} (\infty )= \cC$
and $\varphi(R)=0$.
\end{lem}

\begin{proof}
In lack of a reference we provide a short proof.
Since for every $s\in S$
the curve $\cC_s$ has self-intersection 0 in $\cX_s$,
the cohomology
groups $H^i(\cX_s, \cO_{\cX_s}(\cC_s))$
vanish for $i\ge 1$. Thus for every coherent sheaf
$\cN$ on $S$ the higher direct image sheaves $R^if_*(\cO_X(\cC)
\otimes_{\cO_S}
\cN)$ vanish for $i\ge 1$, see e.g.\
\cite[12.10]{Ha}.
Thus $\cE=f_*(\cO_\cX(\cC))$ is a
locally free sheaf  of rank 2 on $S$, and  forming
$R^0f_*(\cO_\cX(\cC))$ is compatible with
restriction to the fiber, i.e.\ the canonical map
$$
\cE/\fm_s\cE\lto H^0(\cX_s, \cO_{\cX_s}(\cC_s))
$$
is bijective, where $\fm_s$ denotes the ideal sheaf of
the point $s\in S$ (see \cite[12.10 and 3.11]{Ha}).
The inclusion $\cO_\cX\subseteq \cO_\cX(\cC)$ induces
a trivial subbundle $\cO_S$ of $\cE$
(indeed this is true in each fiber). Since the section $R$
is disjoint from $\cC$, the projection
$\cO_\cX\to \cO_R$ extends to a map $\cO_\cX(\cC)\to\cO_R$.
Taking $f_*$ gives a morphism
$\cE\to f_*(\cO_R)\cong \cO_S$ which restricts to the identity
on $\cO_S\subseteq \cE$.
Thus $\cE\cong \cO_S\oplus \cL$ for some line bundle $\cL$ on $S$.
The latter bundle is
trivial due to our assumption that $\Pic (S)=0$. If now $\sigma_0$
and $\sigma_1$ are sections
of $\cE$ which correspond to  the standard basis of
$\cE\cong \cO_S\oplus\cO_S$ then the morphism
$[\sigma_0:\sigma_1]:\cX\to \PP^1$ has the desired properties.
\end{proof}

The following relative version of Castelnouvo's contractibility
criterion is well known\footnote{Cf. e.g., \cite[Theorem
1.3]{KaZa}.}.

\begin{lem}\label{rec.3}
Let $f: \cX \to S$ be a proper smooth family of surfaces and let
$\cC$ be an $S$-flat divisor in $\cX$ such that the fibers
$\cC_s:=f^{-1}(s) \cap\cC$ are smooth rational curves with
self-intersection $-1$ in $\cX_b:=f^{-1}(s)$. Then there exists a
contraction $ \pi : \cX\to \cX'$ of $\cC$, and $\cX'$ is again flat
over $S$.
\end{lem}

\bproof It is sufficient to treat the case where the base $S$ is
affine. In this case there exists an $f$-ample divisor $\cD$ on
$\cX$ which defines an embedding $\cX\hto S\times \PP^N$ for some
$N$. Then the sheaf $\cO_\cX(k\cD-\cC)$, $k\gg 0$, is $f$-semiample
on $\cX$ and provides a desired contraction. \eproof

\blem\label{rec.4} Let $S$ be a scheme with $H^1(S,\cO_S)=0$ and
$\Pic(S)=0$. If $f:\cX\to S$ is a flat morphism with a section
$\sigma:S\to\cX$ such that every fiber is isomorphic to $\PP^1$,
then $\cX$ is $S$-isomorphic to the product $\PP^1\times S$ such
that $\sigma (S)$ corresponds to $\{p\}\times S$ for some point
$p\in \PP^1$. \elem

\bproof We  note first that $R^0f_*(\cO_\cX)\cong\cO_S$ and
$R^1f_*(\cO_\cX)=0$ since the fibers are isomorphic to $\PP^1$.
Using the spectral sequence $H^p(S,R^qf_*(\cO_\cX)) \Rightarrow
H^{p+q}(\cX,\cO_\cX)$ and our assumption $H^1(S,\cO_S)=0$ this
implies that $H^1(\cX,\cO_\cX)=0$. Letting $\Sigma=\sigma(S)$ we
consider the $f$-ample sheaf $\cO_\cX(\Sigma)$. Its direct image
sheaf $\cE=f_*(\cO_\cX(\Sigma))$ is locally free of rank 2 and
$\cX\cong \PP(\cE)$. The sheaf $\cL=\cO_\cX(\Sigma)\otimes
\cO_\Sigma$ is a line bundle on $\Sigma\cong S$ and so is trivial,
since $\Pic(S)=0$ by our assumption. Taking the direct image $f_*$
of the exact sequence
$$
0\to\cO_\cX\to\cO_\cX(\Sigma)
\to\cL\cong \cO_S\to 0
$$
yields an exact sequence
$$
0\to\cO_S\to\cE\to \cO_S\to R^1f_*(\cO_\cX)=0\,.
$$
Thus $\cE$ is an extension of $\cO_S$ by $\cO_S$ and so can be
considered as an element of $\Ext^1_S(\cO_S,\cO_S)$ $\cong
H^1(S,\cO_S)$. Since  by our assumption the latter group vanishes,
this extension splits, i.e.,\ $\cE\cong \cO_S^2$. Hence $\cX\cong
\PP(\cE)= \PP^1 \times  S$, where by our construction $\Sigma$
corresponds to $\{p\}\times S$ for some point $p\in \PP^1$.
\eproof

The following corollary of Lemma \ref{rec.4} is well known;
the proof is immediate.

\bcor\label{rec.51}  Assume that $S$ as in \ref{rec.4} above does
not admit non-constant invertible regular functions. Let $\cC\to
S$ be a flat family of smooth rational curves with a non-empty
$S$-flat subfamily $\cZ\subseteq \cC$ of reduced effective
divisors\footnote{I.e., a disjoint union of images of several
sections $S\to \cC$.}. Then the family $(\cC,\cZ)\to S$ is trivial
i.e., there is an $S$-isomorphism $h:\cC\to \PP^1\times S$ with
$h(\cZ)=\{P_1,\ldots, P_r\}\times S$, where $P_1,\ldots,P_r$ are
points of $\PP^1$.
\ecor

\subsection{Families of completions of a Gizatullin surface}
In this section we  study families of completions of a given
Gizatullin surface $V$. We introduce the notion of a distinguished
extended divisor. In Proposition \ref{propdist} we show that any
deformation family of completions of a Gizatullin surface
over a sufficiently large base is
necessarily trivial provided that the extended divisor is
distinguished and its dual graph stays constant along the
deformation.

\bsit\label{sitfam} We start with the trivial family $f:\cV=V\times
S\to S$, where $S$ is a quasiprojective scheme with $\Pic(S)=0$. We
let $(\bar\cV, \cD)\to S$ be a family of completions of $\cV$ by a
family of standard SNC-divisors $\cD=\bigcup_{i=0}^{n} \cC_i$ over
$S$ with a fixed dual graph. In other words, $\bar\cV\to S$ is a
flat family of complete normal surfaces, $\cD\to S$ is a flat
subfamily of divisors and for every $i$, $f:\cC_i\to  S$ is a flat
family of smooth rational curves which form in every fiber a fixed
standard zigzag (\ref{zigzag}). In particular  for $i=0,\ldots,
n-1$, $\cC_i\cap\cC_{i+1}$ are disjoint sections of $f$.

Since  on the affine part our family is trivial, there is a
simultaneous minimal resolution of singularities $h: \tilde\cV\to
\bar\cV$. This means that $\tilde\cV\to S$ is a smooth family of
complex surfaces, which is fiberwise the minimal resolution of
singularities of $\bar\cV$. Clearly  $h^{-1}(\cV)\cong V'\times S$,
where $V'\to V$ is the minimal resolution.

According to \ref{rec.2} the components $\cC_i$, $i=0,1$,
define morphisms $\Phi_i=
\Phi_{|\cC_i|}:\tilde\cV\to\PP^1$ with
$$
\Phi_0^{-1}(\infty)=\cC_0, \quad
\Phi_1^{-1}(\infty)=\cC_1\quad\mbox{and}
\quad \cC_2\cup\ldots \cup
\cC_n
\subseteq \Phi_0^{-1}(0)\,.
$$
As in the absolute case, we consider the family of divisors
$\cD_{(e)}:= \Phi_0^{-1}(0)$ and the extended divisor
$\cD_\ext:=\cC_0\cup\cC_1\cup \cD_{(e)}$. \esit

It is convenient to introduce the following subgraphs of the extended divisor $D_\ext$ as in (\ref{extended}).

\bsit\label{notsubgraph} For every $1\le i\le n$
we let $D_\ext^{>i}$ denote the union of all
connected components of $D_\ext\ominus C_i$ which do not contain
$C_0$. Similarly we let $D_\ext^{\ge i}$ be the connected
component of $D_\ext\ominus C_{i-1}$ that contains $C_i$.

Obviously,  $D_\ext^{>i}$ is non-empty for every $1\le i\le n-1$,
while $D_\ext^{>n}$ may be empty depending on whether the feather
collection $\{\fF_{nj}\}$ in (\ref{extended}) is empty or not.
\esit

\bdefi\label{dist}
The extended divisor $D_\ext$
will be called {\em distinguished}
if there is no index $i$ with $3\le i\le n$
such that $D_\ext^{>i}$ is non-empty and
contractible.
\edefi

\bprop \label{propdist} Let $V$ be a Gizatullin surface and let
$(\bar\cV,\cD)$  be a family of standard completions of $V$ over
$S=\A^n$ as in \ref{sitfam} with a minimal resolution of
singularities $(\tilde\cV, \cD)$ and extended divisor $\cD_\ext$.
Suppose that at every point $s\in S$ the divisor $\cD_{\ext,s}$ is
distinguished and its dual graph does not depend on $s\in S$. Then the
family $(\bar\cV, \cD)$ is trivial i.e.,  there is an
isomorphism\footnote{Note that this isomorphism is {\em not} the
identity on $V$, in general!} $(\bar\cV, \cD)\cong (\bV, D)\times
\A^m$ compatible with the projection to $\A^m$, where
$\bV=\bar\cV_s$ and $D=\cD_s$ are the fibres over a point $s\in
\A^n$. \eprop

\bproof In the case where $D$ is one of the zigzags $[[0,0]]$ or
$[[0,0,0]]$ the map $\Phi:=\Phi_0\times\Phi_1$ (see \ref{sitfam}) is an isomorphism and the
claim is trivial. Otherwise, since  the dual graph $\cD_{\ext,s}$ at
each point $s\in S$ is the same, we can find a smooth family of
$(-1)$-curves $\cE$ in $\cD_{(e)}$. By Lemma \ref{rec.3} we can
contract $\cE$ simultaneously, which results again in a flat family
of surfaces together with an induced map to
$\PP^1\times\PP^1\times\A^m$. Continuing in this way we get a
sequence of blowdowns \be\la{eq118} \pi:\tilde\cV=\cX_k\to
\cX_{k-1}\to\ldots \to \cX_0=\PP^1\times\PP^1\times \A^m, \ee where
at every step a family of $(-1)$-curves is blown down. Reading this
sequence in the opposite direction, $\tilde \cV$ is obtained from
$\PP^1\times\PP^1\times\A^m$ by a sequence of blowups along sections
say $\Sigma_i\subseteq \cX_i$. Let us show by induction on $i$ that
the family $\cX_i$ is trivial, i.e.\ $S$-isomorphic to $X_i\times
\A^m$ for a suitable blowup $X_i$ of $\PP^1\times\PP^1$. This yields
the desired conclusion, since the triviality of the family
$(\tilde\cV,\cD)$ implies that of $(\bar\cV,\cD)$.

In the case $i=0$ this is evident. If $i=1$
then we can adjust the coordinates in $\PP^1\times\PP^1\times\A^m$
so that the section $\Sigma_1$
is contained in $(0,0)\times \A^m$,
see Lemma \ref{rec.2}. Thus  the first blowup in
(\ref{eq118}) takes place at $(0,0)\times \A^m$ and so  $\cX_1$ is
a trivial family.

Assume by induction that we have an $S$-isomorphism $\cX_i\cong
X_i\times \A^m$ for some blowup $X_i$ of $\PP^1\times\PP^1$.
Let $\cE_j\subseteq \cX_j$ be the exceptional divisor of the $j$th
blowup and $\cE_j^{i}$ its proper transform in $\cX_i$ for $i>j$.
By our assumption the family $\cE_j^i\to S$ is trivial
and $S$-isomorphic to $E_j^i\times\A^m$, where $E_j^i$ is the proper
transform of the $j$th exceptional curve in $X_i$.

If the next blowup is inner with center
$\Sigma_i=\cE_j^i\cap \cE_{j'}^i\approx_S (E_j^i\cap
E_{j'}^i)\times\A^m$, then also $\cX_{i+1}$ is a trivial family.
So assume further that the next blowup is outer with center
$\Sigma_i$ contained in $ \cE_j^i\cong E_j^i\times\A^m$. The
section $\Sigma_i$ is the graph of a map $\sigma_i: \A^m\to E_j^i$
with image contained in
$E_j^i\backslash (D_\ext'\ominus E_j^i)$, where as before
$D_\ext'$ denotes the image of $D_\ext$ in $X_i$. If $E_j^i$ meets
two other components of $D_\ext'$ then $\sigma_i$ maps $\A^m$ to
$\PP^1$ with at least 2 points deleted and so must be constant.
Hence $\cX_{i+1}$ is again a trivial family.

Finally consider the case where $E_j^i$ meets just one other
component of $D_\ext'$. According to Proposition \ref {ext}
all blowups in (\ref{eq118}) are done at the images of the zigzag
$D$.
Thus $E_j^i$ is the image in $X_i$ of some component $C_l$ of $D$.
By our assumption $E_j^i=C_l'$  is an end component of $D_\ext'$,
and so the image $D'$ of $D$ in $X_i$ is a linear chain with end
components $C_0'$ and $C_l'$. Therefore $C_l'$ meets a component
$C_j'$ with $j<l$. Consequently the divisor $D_{\ext}^{>l}$ is
contracted in $X_i$, hence it is contractible.
Since by our assumption $D_{\ext}$ is distinguished this forces
$k=l=2$, and furthermore $i=1$. The latter case was already
treated.
\eproof

In the next Sections 2 and 3 we will show that the condition of
constancy of the dual graph of $\cD_{\ext,s}$ in Proposition
\ref{propdist} is satisfied under the assumptions of Theorem
\ref{MT}. However, in general this condition does not hold as
feathers can jump in families of Gizatullin surfaces. We illustrate
this below by the example of Danilov-Gizatullin surfaces. In Section
2 we will provide a more thorough treatment of this phenomenon.

\bexa
\label{dgsu} Recall that a Danilov-Gizatullin
surface $V=V_{k+1}$ is the complement of a section say
$\sigma$
in a
Hirzebruch surface with self-intersection $\sigma^2=k+1$.
By a theorem of
Danilov-Gizatullin \cite{DaGi}
the isomorphism class of $V_{k+1}$ depends only
on $k$ and not on the choice of $\sigma$ or of the concrete
Hirzebruch surface.  This surface $V_{k+1}$
can be completed by the zigzag
$[[0,0,(-2)_k]]$ with components say $C_0,\ldots, C_{k+1}$.
According to Proposition 5.14 in
\cite{FKZ2}, $V_{k+1}$ admits exactly
$k$ pairwise non-conjugate $\C^*$-actions. In terms of the DPD
presentation
(see \cite{FlZa1} or Section 3 below),
for a fixed $k$ these $\C^*$-surfaces
are given
by the  pairs of
$\Q$-divisors on $C=\A^1$
$$(D_+,D_-)=\left(-\frac{1}{r}[0],\,
-\frac{1}{k+1-r}[1]\right), \qquad r=1,\ldots,k\,.$$ So any other
$\C^*$-action on $V_{k+1}$ is conjugate to one of these.

Given $r\in\{1,\ldots,k\}$ such a $\C^*$-surface $V(r)$ admits an
equivariant standard completion $(\bV(r), D)$ with extended graph
\be\label{DGext} D_\ext(r):\qquad \cou{C_0}{0}\lin\cou{C_1}{0}\lin
\cou{C_2}{-2}\lin\ldots\lin\cou{\quad\qquad C_{r+1}}{-2}
\nlin\cshiftup{\fF_1}{-r}
\vlin{20}\cou{C_{r+2}}{-2}\lin\ldots\lin\cou{\quad\qquad
C_{k+1}}{-2} \nlin\cshiftup{\fF_0}{-1}\qquad\quad \ee where the
curve $C_{r+1}$ is attractive for the extended $\C^*$-action on $\bV(r)$
(cf.\ Section 3). Here the bottom line
corresponds to the boundary zigzag $D$, the feather $\fF_1$ consists
of a single $(-r)$-curve $F_1$ attached to the component $C_{r+1}$
and $\fF_0$ represents a single $(-1)$-curve $F_0$ attached to
$C_{k+1}$. For $r=k$ both feathers $\fF_0,\fF_1$ are attached to the
component $C_{k+1}$. The standard morphism
$\Phi:\bV\to\PP^1\times\PP^1$ is equivariant with respect to a
suitable $\C^*$-action on $\PP^1\times\PP^1$ fixing $(0,0)$.

We note that the extended divisor $D_\ext=D_\ext(r)$ is not
distinguished (see \ref{dist}); indeed,  $D_\ext^{>k+1}=\cF_0$ is
contractible.

Let us construct a family of standard completions of $V_{k+1}$ in
which $D_{\ext,s}$ jumps from one of these extended graphs to
another one, so that $r$ jumps. We restrict for simplicity  to the
case where $r=k$. Blowing down the contractible divisor
$F_0+C_{k+1}$ in $\bV(k)$ we get a new surface $X$ in which
$D_\ext$ is contracted to a chain $[[0,0,(-2)_{k-2}, -1, -k+1]]$
consisting of the images $\bC_0, \ldots, \bC_{k}, \bF_1$ of $C_0,
\ldots, C_{k}, F_1$, respectively. The affine curve
$S:=\bC_{k}\backslash \bC_{k-1}$
is isomorphic to $\A^1$. We let $\cX'$ be the blowup
of the trivial family $X\times S$ along the graph of the embedding
$S\hto  X$ with exceptional curve $\bar\cC_{k+1}$ over $S$.
Finally we let $\cV$ be the blowup of $\cX'$ along a section
$S\hto \bar\cC_{k+1}$ which does not meet the proper transforms of
$\bC_k \times S$ and $\bF_1\times S$, and we denote its
exceptional set by $\cF_0$.

The proper transforms
$$\cC_0,\ldots,\cC_{k+1}, \cF_1
\quad\mbox{of}\quad
\bC_0\times S, \ldots , \bC_k\times S, \bar \cC_{k+1},\bF_1\times S
$$
form together with $\cF_0$ a family of extended divisors $\cD_\ext$ in
$\bar\cV
$, while $\cD=\cC_0\cup\ldots\cup \cC_{k+1}$ is a family of zigzags,
being all of the
same type.

Obviously, the fiber of $(\bar\cV, \cD)$ over the point $s_0$
corresponding to $\bC_k\cap \bF_0$ is $\bV(k)$ with extended
divisor $D_\ext(k)$ while the fibers over the other points $s\in
S\setminus\{s_0\}$ are $\bV(k-1)$ with extended divisor
$D_\ext(k-1)$. Note that all fibers of the  family
$\cV=\bar\cV\backslash \cD\to S$ are isomorphic to $V_{k+1}$ by
the theorem of Danilov-Gizatullin mentioned above. \eexa

\section{Degenerations of singular fibers in families of
$\PP^1$-fibrations}

As we have seen in Example \ref{dgsu}, given a family of standard
completions of a Gizatullin surface $V$ with the same zigzag as in
\ref{sitfam}, the extended divisors $\cD_{\ext,s}$ do not
necessarily have the same dual graph at each point $s\in S$. In this
section we will give a criterion as to when this dual graph stays
constant.

\subsection{Degenerate fibers of a $\PP^1$-fibration}
Let us fix the setup.

\bsit\label{degen.0} Given a surface $V$, we consider a sequence
of blowups \be\label{sbu} \sigma: W=W_m\to W_{m-1}\to \ldots
W_1\to W_0=V \ee with centers in smooth points on $V$ and in its
infinitesimally near points. For $i\ge 1$ we let $E_i\subseteq
W_i$ denote the exceptional $(-1)$-curve of $W_i\to W_{i-1}$. We
consider their proper and total transforms $C_i:=\hat E_i$ and
$\tC_i:=E_i^*$ in $W$, respectively. Clearly the curves $C_i$ (or,
equivalently, the effective cycles $\tC_i$) generate freely the
group $\cy_1(E)$ of 1-cycles supported on the exceptional set $E=
\sum_iC_i$. The intersection form gives a symmetric bilinear
pairing on $\cy_1(E)$. \esit

In the next lemma we describe all cycles  in $\cy_1(E)$
with self-intersection $-1$.

\blem\label{degen.1} (a) $\tC_i.\tC_j=-\delta_{ij}$ for $1\le
i,j\le m$. Moreover $\hF.\tC_i=0$ for every curve $F$ of $V$.

(b) If $C$ is a cycle supported in $E$
with self-intersection $-1$
then $C=\pm \tC_i$ for some $i\ge 1$. In particular
the only effective cycles
with self-intersection $-1$ are the $\tC_i$.

(c) $\tC_i.C_i=-1$ and $\tC_i.C_j\ge 0$ for $i\ne j$.

(d) $\tC_i$ and $\tC_i-C_i$ are orthogonal i.e.,
$\tC_i.(\tC_i-C_i) =0$.
\elem

\bproof To prove (a) we consider the contraction $\pi_i:W\to W_i$,
and we assume that $j\ge i$. If $j>i$ then $\pi_i(\tC_j)$ is a
point and so by the projection formula $\tC_i.\tC_j=\pi_{i*}(C_j)
.E_i=0$. If $i=j$ then with the same argument
$\tC_i.\tC_i=E_i.E_i=-1$. The proof of the second part is similar.

For the proof of (b) we write
$C=\alpha_1\tC_1+\ldots+\alpha_m \tC_m$.
The self-intersection index
$C^2=-\alpha_1^2-\ldots -\alpha_m^2$
is equal to $-1$
if and only if $\alpha_i=\pm 1$ for exactly one $i$
and $\alpha_j=0$ otherwise.

(c) and (d) follow immediately using the projection formula
$\tC_i.C_j=\pi_i^*(E_i).C_j=E_i.\pi_{i*}(C_j)$. \eproof

To study degenerations of extended divisors as introduced in
\ref{unisit}, it is convenient to restrict to the piece
$D_{(e)}=\Phi_0^{-1}(0)$ instead of the full extended divisor
$D_\ext$.

\bsit\label{degen.2}
Letting $\pi:V=U\times\PP^1\to U$,
where $U$ is a neighbourhood of $0\in\A^1_\C$,
we consider
a sequence of blowups as in (\ref{sbu})
with centers on the fiber $F=\{0\}\times \PP^1$
and in infinitesimally near points.
We assume that the full fiber $D_{(e)}=\sigma^{-1}(F)=\hat F+\sum_i
C_i$
has dual graph
\bigskip

\be\label{extended1} D_{(e)}: \qquad\quad \cou{}{\hat F=D_0}
\nlin\xbshiftup{}{\{\fF_{0j}\}} \llin
\ldots\lin\cou{}{D_i}\nlin\xbshiftup{}{\{\fF_{ij}\}}
\llin\ldots\llin \cou{}{D_n} \nlin\xbshiftup{}{\{ \fF_{nj}\}}
\qquad\qquad, \ee where at each curve $D_i$, $0\le i\le n$, a
collection of feathers $\fF_{ij}$ is attached with $1\le j\le
r_i$. Thus each feather $\fF_{ij}$ has dual graph
\be\label{feathergraph}
\co{B_{ij}}\llin\boxo{R_{ij}}\quad\quad =\quad \quad
\co{B_{ij}}\llin\co{R_{ij1}}\llin \ldots \llin \co{R_{ijs_{ij}}}\quad
\quad,
\ee
where the box $R_{ij}$ denotes a linear chain of curves $R_{ijk}$
(possibly empty) connected to the bridge curve $B_{ij}$. We remind
that $R_{ij}$ does not contain a $(-1)$-curve, see Definition \ref
{eboundary.4}. However, unlike in Section 1 we allow that some of
the curves $D_i$ were $(-1)$- curves. This will be convenient in a
later induction argument.

If $D_i$ is one of the curves $C_k$ as considered in \ref{degen.0}
above then we let $\tD_i=\tC_k$. We introduce similarly the
effective cycles $\tB_{ij}$ and $\tR_{ijk}$. Given an irreducible
component  $H$ of one of the feathers $\fF_{ij}$, we call a
component $D_{\mu}$ of the zigzag $D$ a {\it mother component} of
$H$ if $\tilde H.D_{\mu}=1$. \esit

\blem\label{degen.03} \begin{enumerate} [(a)]\item Every component
$H$ of $\fF_{ij}$ has a unique mother component
$D_{\mu}$.
\item  $\tilde H.C=0$ for every component $C$ in
$D_{(e)}$ different from $D_\mu$, $H$ and the neighbor of $H$ in (\ref{feathergraph}) to the right.
\end{enumerate} \elem

\bproof
Let $H$ be the curve $C_k$ considered in
\ref{degen.0} so that
$\tilde H=\pi_k^*(E_k)$, where $E_k$ is the
exceptional $(-1)$-curve
created in the blowup $W_k\to W_{k-1}$ and
$\pi_k:W\to W_k$ is the
contraction as in the proof of \ref{degen.1}.
Since $H=\hat E_k$ does not separate the
zigzag, the center of
$\sigma_k:W_k\to W_{k-1}$ cannot be
a double point of
the zigzag $\pi_k(D)$. Thus $E_k$ meets a unique component of $\pi_k(D)$. In view of the projection formula $\tilde H.D_i=E_k.\pi_k(D_i)$  this proves (a).

To deduce (b), assume that $C\subseteq D_{(e)}$ is a component
different from $D_\mu,H$ and satisfying  $\tilde H.C\neq 0$. Again
by the projection formula $\tilde H.C=E_k.\pi_k(C)\ne 0$. Since
$E_k$ is an at most linear vertex of the dual graph of
$\pi_k(D_{(e)})$, this is only possible if $E_k.\pi_k(C)=1$. As
observed before, $H$ does not separate the zigzag $D$, and so $C$
being different from $D_\mu$ must belong to the feather
$\fF_{ij}$. Since $\tR_{ij}$ contains no $(-1)$-curves, the
projection formula forces $C$ to be the neighbor on the right in
(\ref{feathergraph}) as claimed in (b).
\eproof

\bexa\label{degen.4} To illustrate these notions let us consider the
graph \vspace*{0.5truecm}
$$
D_{(e)}:\qquad\quad \cou{\hF=D_0}{w_0}\lin \ldots \lin
\cou{D_j}{w_j}\lin \cou{D_{j+1}}{-2}\lin \ldots
\lin\cou{D_{l-1}}{-2} \lin\cou{\quad\quad D_l}{-2} \nlin
\cshiftup{B}{\qquad\quad j-l-1}
\lin\!\!\lin\cou{D_{l+1}}{w_{l+1}}\lin \ldots\quad,
$$
 where $D_{(e)}^{\ge l+1}$ is contractible
to a smooth point on $D_l$. It is easily verified
that the mother component
of $B$ is $D_j$.
\eexa

In the next proposition we  collect some important properties of
mother components. For a graph $D_{(e)}$ as in (\ref{extended1}),
similarly as before, $D_{(e)}^{>i}$ denotes the union of all
connected components of $D_{(e)}\ominus D_i$ not containing $D_0$,
while $D_{(e)}^{\ge i}$ stands for the connected component of
$D_{(e)}\ominus D_{i-1}$ containing $D_i$.

\bprop\label{degen.3}
\begin{enumerate}[(a)]
\item Let $D_{\mu}$, $\mu=\mu(i,j)$, be the mother component of
$B_ {ij}$. Then $\mu\le i$ and $B_{ij}^2\ge \mu-i-1$. Furthermore,
$B_{ij}^2\le -1$ and  $B_{ij}^2=-1$ if and only if $\mu=i$.
\item If $\mu<i$  then the
divisors \be\label{restcur} D_{(e)}^{\ge \mu+1}\ominus\fF_{ij},\quad
D_{(e)}^{>i}\ominus\fF_{ij}\quad\mbox{and}\quad
\fF_{i'j'}\quad\mbox{with}\quad \mu<i'< i \ee are all contractible
inside $D_{(e)}$.

\item Assume that $B_{i'j'}$ (where
$(i',j')\neq (i,j)$) is a further bridge curve with mother
component $D_{\mu'}$. If $\mu< i$ and $\mu'< i'$ then the
intervals $[\mu+1,i]$ and $[\mu'+1, i']$ are disjoint.
\end{enumerate}
\eprop

\bproof
(a) The piece of the zigzag $D$ between
$D_{\mu+1}$ and $D_i$
separates the mother component $D_{\mu}$
from $B_{ij}=\hat E_k$.
Hence $\pi_k$ contracts this piece. Furthermore,
$D_{\mu}$
separates it from $D_0=\hat F$ within $D$,
therefore $\mu\le i$.
Moreover at most $i-\mu$ blowups are done near
the images of $B_{ij}$, so $B_{ij}^2\ge \mu-i-1$.

The remaining assertions of (a) are
easy and so we leave the proof to the reader.

(b) To show that $D_{(e)}^{\ge \mu+1}\ominus\fF_{ij}$ is
contractible, it suffices to verify that $D_{(e)}^{\ge \mu+1}\ominus
R_{ij}$ supports the total preimage $\tB_{ij}$, since then it
contracts to $E_k=\pi_k(B_{ij})$ under $\pi_k$. As before
$E_k=\pi_k(B_{ij})$ represents an at most linear vertex of the dual
graph of $\pi_k (D_{(e)})$, where one neighbor is $D_\mu$ and the
other one (if existent) is the neighbor of $B_{ij}$ in $\fF_{ij}$ to
the right in (\ref{feathergraph}). Moreover all components in
$D_{(e)}^{\ge \mu+1}\ominus\fF_{ij}$ appear under further blowups
with center at $\pi_k(D_\mu)\cap E_k$ and its infinitesimally near
points. Hence the assertion follows.

If $D_{(e)}^{>i}\ominus\fF_{ij}$ were not contractible,
then contracting successively all $(-1)$-curves in
$D_{(e)}^{\ge \mu+1} \ominus\fF_{ij}$ the vertex $D_{i}$ of
(\ref{extended1}) would remain a branching point, which is
impossible.

Similarly, if a feather $\fF_{i'j'}$ with $\mu<i'< i$ were not
contractible, then contracting successively all $(-1)$-curves in
$D_{(e)}^{\ge \mu+1} \ominus\fF_{ij}$, the vertex $D_{i'}$ of
(\ref{extended1}) would remain a branching point, which is
impossible.

To show (c) we let $B_{ij}=\hat E_k$ and $B_{i'j'}=\hat E_{k'}$.
We may assume that $k'<k$ so that $\pi_k$ does not contract
$B_{i'j'}$. As the divisor $D_{(e)}^{\ge \mu+1}\ominus\fF_{ij}$
is contracted under $\pi_k$
this implies that $i'\le\mu$. Hence
$$
\mu'< i'\le \mu <i,
$$
proving (c).
\eproof

\subsection{Families of rational surfaces: the specialization map}
Let us recall the notion of specialization and generalization
map for smooth proper families.

\bsit\label{spemap} We consider a proper smooth holomorphic map
$\pi: \cX\to\Di$, where $\cX$ is a connected complex manifold and
$\Di$ stands for the unit disc in $\C$ with center $0\in\Di$. By
Ehresmann's theorem for any point $s\in \Di$ there is a
$\Di$-diffeomorphism $\cX\cong \cX_s\times \Di$, where
$\cX_s=\pi^{-1} (s)$. Hence the embedding $\cX_s\hookrightarrow
\cX$ induces an isomorphism in cohomology
$H^*(\cX)\stackrel{\cong}{\longrightarrow} H^*(\cX_s)$. Composing
the isomorphisms
$$H^*(\cX_s)\stackrel{\cong}{\longrightarrow}
H^*(\cX)\stackrel{\cong}{\longrightarrow} H^*(\cX_0)$$
we obtain a {\it specialization map}
$\sigma:H^*(\cX_s)\stackrel{\cong}{\longrightarrow} H^*(\cX_0)$;
its inverse is called
a {\it generalization map}.
\esit

\bsit\label{rasu} From now on  we assume that the fibers $\cX_s$
are complete rational surfaces.
Then
$$\NS (\cX_s)=\Pic (\cX_s)\cong H^2  (\cX_s;\Z)\,,$$
where $\NS (\cX_s)=\Div (\cX_s)/\sim$ is the Neron-Severi group of
algebraic 1-cycles modulo numerical  equivalence. From the exact
 sequence
$$0=H^1 (\cX,\cO_\cX) \to H^1 (\cX,\cO_\cX^\times)\cong \Pic (\cX)
\to H^2 (\cX,\Z)\to H^2 (\cX,\cO_\cX)=0$$ induced by the
exponential sequence we obtain an isomorphism $$\NS (\cX)=\Pic
(\cX)\cong H^2  (\cX;\Z)\,,$$ which commutes with restrictions to
the fibers that is, with the isomorphisms \be\label{piiso} \Pic
(\cX) \stackrel{\cong}{\longrightarrow} \Pic (\cX_s)
\qquad\mbox{and}\qquad H^2  (\cX;\Z)
\stackrel{\cong}{\longrightarrow} H^2   (\cX_s;\Z)\ee induced by
the embeddings $\cX_s\hookrightarrow \cX$. Composing the above
isomorphisms leads to
$$\sigma: \NS (\cX_s) \stackrel{\cong}{\longrightarrow}
\NS (\cX_0)\,$$  also called a {\it specialization map}. Clearly
$\sigma$ is an isometry with respect to the intersection forms.
\esit

\blem \label{effcon} For a general  point $s\in\Di$, the
specialization map $\sigma $ sends the effective cone in $\NS
(\cX_s)\otimes \Q$ into the effective cone in $\NS (\cX_0)\otimes
\Q$. \elem

\bproof
For an invertible sheaf $\cL\in\Pic (\cX)$,
its direct image
$R^1\pi_*(\cL)$ is a coherent sheaf on $\Di$,
with a torsion located on a discrete set, say,
$A(\cL)\subseteq \Di$. Since
$\Pic_0 (\cX)=0$
the set $A=\bigcup_{\cL\in\Pic (\cX)} A(\cL)$
is at most countable.

Picking now a point $s\in\Di\setminus A$, for an effective 1-cycle
$C$ on $\cX_s$ we consider the corresponding invertible sheaf
$\cL_s=\cO_{\cX_s} (C)$. By virtue of (\ref{piiso}) there exists an
invertible sheaf  $\cL\in\Pic (\cX)$ such that $\cL\vert
\cX_s=\cL_s$. Since $R^1\pi_*(\cL)$ has no torsion at $s$, the
restriction map $$H^0(\cX,\cL) \to H^0 (\cX_s,\cL_s)$$ is a
surjection, and so the sections of the sheaf $\cL_s$ can be lifted
to sections of $\cL$. In particular $\cL=\cO_{\cX} (\cC)$ for some
effective 1-cycle $\cC$ on $\cX$ with $\cC|\cX_s = C$. Hence also
$\sigma (C)=\cC\vert \cX_0$ is effective. This yields the lemma.
\eproof

\subsection{Formal specialization map and jumping feathers}
In this section we study possible degenerations of families of
extended divisors. We  recall first the geometric setup of Section
1.4.

\bsit\label{degen.51} Let $V=X\setminus D$ be a Gizatullin surface
with a boundary zigzag $D$. As in Section 1.4 we consider families
of standard completions $(\tilde\cV,\cD_s)$, $s\in S $, of a
minimal resolution of singularities $V'\to V$  with a
corresponding family of extended divisors
$(\cD_\ext)_s=(\cC_0)_s+(\cC_1)_s+(\cD_{(e)})_s$. We are
interested in degenerations in such families. More precisely, each
divisor $(\cD_{(e)})_s$ has a dual graph as in (\ref{extended1}),
{\em however this graph may depend on} $s\in S$. If
$\fF_{ij}(s)=\cB_{ij}(s)+\cR_{ij}(s)$ denotes the feathers at the
point $s$ then clearly the part $\cR_s=\sum \cR_{ij}(s)$ must be
constant being the exceptional set of the resolution of
singularities of $V$. Similarly the  dual graph of the boundary
zigzag $\cD_s\cong D$ stays constant.

Assuming that $S$ is a smooth curve,
for a general
point $s\in S$ the  specialization map
$$\sigma: \NS (\tilde\cV_s) \stackrel{\cong}{\longrightarrow}
\NS (\tilde\cV_{s_0})\,$$ restricts to an isomorphism
$$\sigma: \NS ((\cD_{(e)})_s)=\cy_1 ((\cD_{(e)})_s )
\stackrel{\cong}{\longrightarrow} \NS ((\cD_{(e)})_{s_0})=\cy_1
((\cD_{(e)})_{s_0})\,$$ of the corresponding cycle spaces
compatible with the intersection forms. In what follows we study
this map $\sigma$ on a  formal level. \esit

\bsit\label{degen.6} Let us consider two modifications $\pi:W\to
V$ and $\pi':W'\to V$ as in \ref{degen.2} above, with the same
number $m$ of blowups. Moreover assume that on $W,W'$ we have
 decompositions $$D_{(e)}=D+ \sum \fF_{ij}\qquad\mbox{and}\qquad
 D_{(e)}'=D'+ \sum \fF'_{ij}$$ as in \ref{degen.2}
with the same number $n$ of curves $D_i,D'_i$ and with feathers
$$\fF_{ij} =\;\; \co {B_{ij}} \lin\boxo{R_{ij}}
\qquad\mbox{and}\qquad \fF'_{ij} =\;\; \co {B'_{ij}}
\lin\boxo{R'_{ij}}\quad,$$ respectively. We let $G=\cy_1
(D_{(e)})$ and $G'=\cy_1 (D_{(e)}')$ be their groups of 1-cycles
with generators $(C_i)$ and $(C_i')$ or, equivalently, $(\tC_i)$
and $(\tC_i')$, respectively. Suppose that we are given an
isomorphism
$$
\delta:G\to G'
$$
with the following properties:
\begin{enumerate}[(i)]
\item $\delta$ respects the intersection forms.
\item $\delta$
transforms effective cycles into effective cycles. \item
$\delta(D_i)=D_i'$ for all $i$. \item $\delta(R_{ijk})=
R'_{i'j'k}$ for some $i', j'$, where $R_{ijk} $, $ R'_ {i'j'k}$
are the components of $R_{ij}$, $ R'_{i'j'}$, respectively,
ordered as in \ref{degen.2}.
\end{enumerate}
We then call $\delta$ a {\it formal specialization map},
and $\delta^{-1}$ a
{\it formal generalization map}.

It is clear from the discussion in \ref{degen.51} that any
specialization map arising from  a degeneration in a family of
completions/resolutions of a Gizatullin surface is also a formal
specialization map. Indeed (i) and (ii) follow from the construction
in view of Lemma \ref{effcon}, (iii) follows immediately by the
triviality of the family $\cD\to S$, and (iv) holds due to the
constancy of singularities in the open part $\cV_s\cong V$. \esit
{\em We assume in the sequel  that $\delta$ is a formal
specialization map.}

The structure of $\delta$  can be understood on the level
of the generators $\tD_i$, $\tR_{ijk}$
and $\tB_{ij}$ of $G=\cy_1 ((\cD_{(e)})_s)$.
These generators form an orthogonal basis of $G$
(see Lemma \ref{degen.1}(a)).
The same is true for their images
in $G'$. So according
to Lemma \ref{degen.1}(b)
\be\label{-1curves}
\{\delta(\tD_i),\delta(\tB_{ij}),\delta(\tR_{ijk})\}=
\{\tD'_i,\tB'_{ij},\tR'_{ijk}\}\,.
\ee

\bprop\label{degen.7} With the assumptions as before the following
hold.
\begin{enumerate} [(a)]
\item $\delta(\tD_i)=\tD_i$ and
$\delta(\tR_{ijk})=\tR'_{i'j'k}$\,; \item
$\delta(\tB_{ij})=\tB'_{i'j'}$\,; \item $\delta$ respects the
mother components, i.e.\ if $D_{\mu}$ is the mother component
of $B_{ij}$ then $D_{\mu }'$ is the mother component of
$B'_{i'j'}$. \item Every feather $\fF_{ij}=B_{ij}+R_{ij}$ either
stays fixed or jumps to the right under $\delta$, i.e.\
$\delta(\tB_{ij})=\tB'_{i'j'}$ and $\delta(R_{ij})=R'_{i'j'}$ with
$i'\ge i$.
\end{enumerate}
\eprop

\bproof To deduce (a) we note  that by \ref{degen.1}(b)
$\delta(\tD_i)= \tC'$ for some irreducible component $C'$ of
$D'_{(e)}$. Using properties (i) and (iii) of $\delta$
$$
\tC'.D_i'= \delta(\tD_i).D'_i=\tD_i.\delta^{-1}(D_i')=
\tD_i.D_i=-1.
$$
Using Lemma \ref{degen.1}(c) this implies that $C'=D'_i$. With the
same argument
it follows that $\delta(\tR_{ijk})=\tR'_{i'j'k}$. Clearly (b) is a
consequence
of (a) and (\ref{-1curves}).

(c) follows from the equation
$$
\tB_{ij}.D_ \alpha= \delta(\tB_{ij}).\delta(D_\alpha)=\tB_{i'j'}'.
D'_\alpha
$$
and the characterization of mother components given in Lemma \ref
{degen.03}.

(d) By property (ii) $\delta$ sends the effective cone of $G=\cy_1
(D_{(e)})$ into the effective cone of $G'=\cy_1 (D_{(e)}')$.
Moreover, $B_{ij}$, $B'_{i'j'}$ appear in the cycles $ \tB_{ij}$
and $\delta(\tB_{ij})=\tB'_{i'j'}$, respectively, with coefficient
1. Hence $\delta (B_{ij})=B_{i'j'}'+\Delta$ with an effective
divisor $\Delta=\Delta(i,j)$ which does not contain $B_ {i'j'}'$
so that
$$
\Delta=\sum_{p\ge 1}\alpha_p D_p' + \sum_{(p,q)\neq
(i',j')}\alpha_{pq} B_{pq}'+ \sum_{p,q,r}\alpha_{pqr} R_{pqr}' \,,
\quad \mbox{where}\quad \alpha_{p},\,\alpha_{pq},\, \alpha_{pqr}
\ge 0\,.$$ Suppose that $B_{ij}$ jumps indeed, i.e.\ $i\ne i'$.
Then $D_{i'}.B_{ij}=0$, hence
$$
0=D'_{i'}.\delta (B_{ij})=D'_{i'}.B'_{i'j'}+D'_{i'}.\Delta=1+D'_{i'}.
\Delta.
$$
Thus $D'_{i'}.\Delta=-1$ and so $\alpha_{i'} >0$. It follows that
$K:=\{p\, : \,\alpha_p >0\}$ contains $i'$. It is easily seen that
$0\notin K$. We choose $p\in\{0,\ldots,n\}\setminus K$ so that at
least one of $p\pm 1$ is in $K$. Since  $\delta(D_i)=D_i'$
$\forall i$
and $\delta$ preserves the intersection form, we have
\be\label{nazv} D_{p}. B_{ij}=\delta(D_{p}).
\delta(B_{ij})=D_{p}'. B_{i'j'}'+D_{p}'.\Delta\ge
D_{p}'.\Delta>0\,. \ee Hence $D_{p}. B_{ij}=1$ and so $p=i$.
Consequently $K=[i+1,\ldots,n]$  (indeed, $0\notin K$). Since
$i'\in K$ this proves (d). \eproof

\subsection{Rigidity}
In Theorem \ref{maincrit} below we give a criterion for the dual
graph of the extended divisor $D_\ext$ to stay constant under any
specialization or generalization. We use the following
terminology.

\bdefi\label{rig} We say that the divisor $D_{(e)}$ as in
\ref{degen.2} is {\em stable under specialization} if for any
specialization map $\delta: G=\cy_1(D_{(e)})\to
G'=\cy_1(D'_{(e)})$ as in \ref{degen.6} we have
$\delta(B_{ij})=B'_{ij}$ with a suitable numbering of $B'_{i1},
\ldots, B'_{i, r_i}$. This means that no feather jumps to the
right in (\ref{extended1}).

Similarly, a divisor $D_{(e)}$ is said to be {\em stable under
generalization} if for any generalization map\footnote{i.e.,
$\delta=\gamma^{-1}$ is a specialization map as in \ref{degen.6}.}
$\gamma: G=\cy_1(D_{(e)})\to G'=\cy_1(D'_{(e)})$ we have
$\gamma(B_{ij})=B'_{ij}$ with a suitable numbering of $B'_{i1},
\ldots, B'_{i, r_i}$. Therefore no feather jumps to the left in
(\ref{extended1}).

Finally, a divisor $D_{(e)}$, which is stable under both
specialization and generalization, is said to be {\em rigid}. This
terminology can be equally applied to the extended divisor
$D_\ext=C_0+C_1+D_{(e)}$. \edefi

We have the following fact.

\bprop\label{corgen}
$D_{(e)}$ is stable under generalization if and only if
$B^2_{ij}=-1$ for all bridge curves $B_{ij}$. \eprop

\bproof  By Proposition \ref{degen.7}(d) $B_{ij}$ can only jump to
the left under generalization so that
$\gamma(\tB_{ij})=\tB'_{i'j'}$ with $i'\le i$. Assume that
$B^2_{ij}=-1$. Then $D_i$ is the mother component of $B_{ij}$,
see Proposition \ref{degen.3}(a). By virtue of
Proposition \ref{degen.7}(c) $D'_i$ is the mother
component of $B_{i'j'}'$. Using again Proposition \ref{degen.3}(a)
$i\le i'$, hence $i=i'$ and the feather $\fF_{ij}$ stays
fixed, as required.

To show the converse we assume that $B^2_{ij}\le -2$. By
Proposition \ref{degen.3}(a) then $\mu<i$, where
$D_\mu$ is the mother component of $\fF_{ij}$. Using Proposition
\ref{degen.3}(b) the divisor $P:=D_{(e)}^{\ge \mu +1} \ominus \fF_
{ij}$ is contractible. Let $F$ be the contracted divisor
$D_{(e)}/P$. Then the image of $P$ is the intersection point of
the images, say $ \bB_{ij}$ and $\bD_\mu$ of $B_{ij}$ and $D_\mu$
in $F$. Moreover $\bB_{ij}^2=-1$.

Let now $q$ be a point on $\bD_\mu$ different from this
intersection point. Rebuilding $P$ at this point $q$ yields a new
divisor, say, $D'_{(e)}$. We claim that this procedure provides a
non-trivial generalization map $\gamma: \cy_1(D_{(e)})\to
\cy_1(D'_{(e)})$ or, equivalently, a non-trivial specialization
map $\delta: \cy_1(D'_{(e)})\to \cy_1(D_{(e)})$.

Obviously the curves in $D_{(e)}$ and $D'_{(e)}$ are in
$1-1$-correspondence. Let for a curve $C$ in $D_{(e)}$, $C'$
denote the corresponding curve in $D'_{(e)}$. We define $\delta$
by $\delta(C')=C$ for $C'\ne B'_{ij}$ and
$\delta(B_{ij}'):=\tB_{ij}$. Since
$$
\tB_{ij}. C=0\quad\mbox{for }C\ne B_{ij}, D_\mu, \qquad
\tB_{ij}^2=-1\quad\mbox{and}\quad \tB_{ij}.D_\mu=1\,,
$$
$\delta$ is an isometry. Since it maps effective cycles into
effective cycles, $\delta$ is a specialization map, as required.
\eproof

\bprop\label{corspec} Assume that a feather $\fF_{ij}$ in
(\ref{extended1}) jumps to
$\fF'_ {i'j'}$ under a specialization $\delta: D_{(e)}\to
D'_{(e)}$. If $i'> i$ then the following divisors are either empty
or contractible:

(a) $\fF_{kl}$ with $i<k< i'$;

(b) $D_{(e)}^{>i'}$ and $D_{(e)}^{\prime\; >i'}\ominus \fF'_ {i'j'}$;

(c) $D_{(e)}^{\ge i+1}$ and $D_{(e)}^{\prime\;\ge i+1} \ominus \fF'_
{i'j'}$.
\eprop

\bproof (a) follows from Proposition \ref{degen.3}(b).
Indeed, if $D_\mu$ is the mother component of $B_{ij}$ the by
Proposition \ref{degen.7}(c) $D'_\mu$ is the mother component of $B'_{i'j'}$, so
$\mu\le i<k<i'$. Similarly by the same
Proposition \ref{degen.3}(b),
$D_{(e)}^{\prime \;>i'}\ominus \fF'_ {i'j'}$ is either empty or
contractible, as stated in (b).

Now (b) and (c) can be shown by induction on the number of
irreducible components of $D_{(e)}$. Let (b)$_m$ and (c)$_m$ be
the corresponding statements for divisors with $m$ components. We
show below that
\smallskip

(i) (b)$_{m-1}$,(c)$_{m-1}\Rightarrow$ (b)$_m$, and

(ii) (b)$_{m}$,(c)$_{m-1}\Rightarrow$ (c)$_m$.
\smallskip

\noindent To deduce (i) and (ii) we use the following claim.

{\it Claim 1. Suppose that the divisor $D_{(e)}^{> i'}$ is
non-empty. Then there exist $(-1)$-curves $C$ in $D_{(e)}^{> i'}$
and $C'$ in $D_{(e)}^{\prime\; > i'}$ with $\delta(C)=C'$, which are
contractible in  $D_{(e)}$ and $D'_{(e)}$, respectively.}

The contractibility of $D_{(e)}^{> i'}$, and then also (i) and (b),
follow from this claim by induction on $m$. Indeed, contracting
$C,C'$ in $D_{(e)}$, $D'_{(e)}$, respectively, leads to new
divisors, say, $ D_{(e)}^\vee$ and $D_{(e)}^{\prime\,\vee}$, where
$D_{(e)}^{\prime\,\vee}$ is a specialization of $D_{(e)}^\vee$. By
virtue of Proposition \ref{degen.3}(a) the feather $\fF'_ {i'j'}$ is
minimal. Hence $\fF_{ij}=\fF_{ij}^\vee$ and $\fF'_ {i'j'}={\fF'}_
{i'j'}^\vee$ are not affected by these contractions and again
$\fF_{ij}^\vee$ jumps to ${\fF'}_ {i'j'}^\vee$.

\noindent {\it Proof of Claim 1.} Assume first that $i'<n$. The
divisor $D_{(e)}^{\prime\; \ge i'+1}$ is then non-empty and
contractible. Hence it contains a $(-1)$-curve $C'$ representing an
 at most linear vertex of the dual graph of ${D_{(e)}^\prime}$.
This curve $C'$ can be either $D'_{k'}$ or a bridge $B'_{k'l'}$,
where $k'\ge i'+1$.

In the latter case we let $\tB'_{k'l'}=\delta(\tB_{kl})$. Since
$(B'_{k'l'})^2=-1$, by Propositions \ref{degen.3}(a) and
\ref{degen.7}(c) $D_{k'}$ and $D'_{k'}$ are the mother components
of $B_{kl}$ and $B'_{k'l'}$, respectively. Hence $k\ge k'$. Since
 under specialization
a feather can only jump to the right (see Proposition
\ref{degen.7}(d)), we have $k=k'$. Therefore again by Proposition
\ref{degen.3}(a), $B_{kl}^2=-1$ and the curves $C=B_{kl},
\,C'=B'_{k'l'}$ are as desired. Indeed, in view of Lemma
\ref{degen.1}(b), $B_{kl}=\tB_{kl}$, $B'_{k'l'}=\tB'_{k'l'}$ and
so  by Proposition \ref{degen.7}(b), $\delta(B_{kl})= B'_{k'l'}$.

In the former case $C=D_{k'}$ is again a $(-1)$-curve, since
$\delta$ respects the intersection forms  and  $C'=\delta(C)$. If
$D_{k'}$ is at most linear vertex of the dual graph then the
curves $C,C'$ are as desired. Otherwise $D_{k'}$ is a branch point
of the dual graph  while $D'_{k'}$ is not. So there is a feather
$\fF_{k'l}$ at $D_{k'}$ which jumps to the right under $\delta$.
Thus $k'<n$, and we can repeat the consideration using induction
on $k'$.

Suppose further that $i'=n$. Since by our assumption $D_{(e)}^{> n}$
is non-empty, there is a non-empty feather, say, $\fF_{nl}$ at
$D_n$. This feather stays fixed under $\delta$ i.e., $\delta
(\tB_{nl})=\tB'_{nl'}$. Moreover, since $\mu<n$, by virtue of
Proposition \ref{degen.3}(c) $D_n$ and  $D_n'$ are the mother
components of $\fF_{nl}, \fF'_{nl'}$, respectively. Similarly as
above, this implies that $C=B_{nl}$, $C'=B'_{nl'}$ are $(-1)$-curves
with $\delta(C)=C'$, as desired. This proves the claim.
\smallskip

The proof of (ii) proceeds in a similar way. Because of (b)$_m$ we
may assume that $D_{(e)}^{>i'}$ and
$D_{(e)}^{\prime\; >i'}\ominus \fF'_{i'j'}$ are empty since
otherwise we can contract them inside $D_{(e)}$, $D'_{(e)}$,
respectively, and use induction on $m$ as before. Similarly due to
(a) we may suppose that both $D_{(e)}$ and $D'_{(e)}$ have no
feathers at components $D_k$, $D'_k$ with $i<k<i'$.

Now the induction step can be done due to the following
\smallskip

{\it Claim 2. Under the assumption as above there are $(-1)$-curves $C$ in $D_{(e)}^{\ge i+1}$ and $C'$ in $D_{(e)}^{\prime\; \ge i+1}\ominus
\fF'_{i'j'}$,
with $\delta (C)=C'$, which are contractible in $D_{(e)}$, $D'_{(e)}$,
respectively. }
\smallskip

\noindent
{\it Proof of Claim 2.} These divisors in our case consist of the
linear strings $[D_{i+1},\ldots,D_{i'}]$ and
$[D'_{i+1},\ldots,D'_{i'}]$, respectively. It is enough to show
that there is a $(-1)$-curve in one of these linear strings, and
then similarly as above there is also the second one.

Let as before $D_\mu$ be the mother component of the bridge curve
$B_{ij}$. Then $D'_\mu$ is the mother component of $B'_{i'j'}$. If $
\mu=i (<i')$ then by Proposition \ref{degen.3}(b) the non-empty
divisor $D_{(e)}^{\prime \ge i +1}\ominus
\fF'_{i'j'}=[D'_{i+1},\ldots,D'_{i'}]$ is contractible and so the
result follows. If $\mu <i$ then again by Proposition
\ref{degen.3}(b) the divisor $D_{(e)}^{>i}\ominus \fF_{ij}$ is
contractible, and also its connected component $D_{(e)}^{\ge
i+1}=[D_{i+1},\ldots,D_{i'}]$ is. Hence again we are done, and so
the proof is completed. \eproof

The following fact is in a sense a converse to Proposition
\ref{corspec}.

\bprop\label{convspec} Suppose that, for two indices $i, i'$ with
$0\le i<i'\le n$, each one of the following divisors is either
empty or contractible:

(a) the feathers $\fF_{kl}$ with $i<k< i'$;

(b) the divisor $D_{(e)}^{>i'}$;

(c) the divisor $D_{(e)}^{\ge i+1}$.

\noindent Then any feather $\fF_{ij}$ jumps to a feather
$\fF'_{i'j'}$ under a suitable specialization. \eprop

\bproof
The proof is similar to that of Proposition
\ref{corgen}. Contracting first $D_{(e)}^{>i'}$ and then the
remaining part, say, $P$ of $D_{(e)}^{\ge i+1}$, we rebuild $P$
blowing up at the intersection point of $D_i$ and $\fF_{ij}$ and
its infinitesimally near points. After that we rebuild
$D_{(e)}^{>i'}$ at points of $D_{i'}$ different from the
intersection point with the new feather $\fF'_{i'j'}$. We leave
the details to the reader.
\eproof

Now we are ready to formulate our main rigidity criterion. This
enables us in the next section to check rigidity for Gizatullin
$\C^*$-surfaces satisfying one of the conditions ($\alpha_+$),
($\alpha_*$) or ($\beta_+$), ($\beta_*$) of Theorem \ref{MT}.

Similarly as in \ref{dist} we call a divisor $D_{(e)}$ {\it
distinguished} if there is no index $i$ with $1\le i\le n$ such
that $D_{(e)}^{>i}$ is non-empty and contractible.

\bthm\label{maincrit} A distinguished divisor
$D_{(e)}$ is rigid provided that
all its bridges $B_{ij}$ are $(-1)$-curves and one of the
following conditions is satisfied.
\bnum[(i)] \item $D_{(e)}^{>
n}\ne\emptyset$.
\item If for some $i$, $0\le i< n$, the feather
collection $\{\fF_{ij} \}$ is non-empty then the divisor
$D_{(e)}^{\ge i+1}$ is not contractible. \enum \ethm

\bproof By Proposition \ref{corgen} $D_{(e)}$ is stable under
generalization. Suppose on the contrary that a feather $\fF_{ij}$
jumps to $\fF'_{i'j'}$ under a specialization, where $i<i'\le n$. By
Proposition \ref {corspec}(c) $D_{(e)}^{\ge i+1}$ is contractible
and so (ii) is violated. Similarly, by Proposition \ref {corspec}(b)
$D_{(e)}^{> i'}$ is contractible. Since  $D_{(e)}$ is distinguished
and $i'+2\ge 3$, this is only possible if  $i'=n$ and $D_{(e)}^{>
n}=\emptyset$. Thus (i) is violated as well, proving the theorem.
\eproof

We finish this section with several examples of rigid or non-rigid
divisors.

\bexas\label{exa22} 1. Consider the  Gizatullin $\C^*$-surface $V$
defined by the following pair of $\Q$-divisors on $\A^1$ (see
Section 3.1):
$$(D_+,D_-)=\left(\frac{1}{n}[0]-[1],
\,-\frac{1}{n}[0]\right)\,.$$ According to Proposition
\ref{eboundary.9} below its standard completion has degenerate fiber
with dual graph \vskip3truemm
$$ D_{(e)}\,:\qquad
\cou{D_0}{-n}\vlin{17}\cou{\!\!\!\!\!\!\!\!\!\!\!\! D_1}{-2}
\nlin\cshiftup{B_1}{-1}\lin\cou{D_2}{-2}\lin\ldots\lin
\cou{D_{n-1}}{-2}\lin\cou{D_n}{-2}\quad.
$$

\noindent Using
Propositions \ref{corgen} and \ref{corspec} the divisor $D_{(e)}$
is rigid i.e., stable under specialization or generalization.

\smallskip

2. Let us revisit the standard completion of a Danilov-Gizatullin
surface $V=V_{n}$ with $n=k+1\ge 3$ (see \ref{dgsu}), which has
extended divisor (\ref{DGext}).
The feather $\fF_1$ has mother component $C_2$. By Proposition
\ref{corgen} it can jump to $C_2$ under a suitable generalization,
but also to any other component $C_i$, $i\ge 2$,  using Proposition
\ref {convspec}.

3. Let $D_{(e)}$ be the divisor

\vskip3truemm
$$D_{(e)}\,:\qquad\cou{D_0}{-n}\lin\!\!\lin
\cou{\!\!\!\!\!\!\!\!\!\!\!\! D_1}{-2}
\nlin\cshiftup{B_1}{-1}\lin \cou{D_2}{-2}
\lin\ldots\lin\cou{}{-2}\nlin\cshiftup{B_2}{-2}
\quad D_{n-1}\quad\quad.
$$

\no Again this is the dual graph of the degenerate fiber in a
standard completion of a Gizatullin $\C^*$-surface
$V=\Spec\,\C[t][D_+,D_-]$ with
$$(D_+,D_-)=\left(\frac{1}{n}[0]-[1],
\,-\frac{1}{n-1}[0]\right)\,$$ (cf.\ Example 1 above). According
to Proposition \ref{corspec}, $D_{(e)}$ does not admit a
nontrivial specialization. The mother component of $B_2$ is $D_0$.
Hence by Proposition \ref{corgen}, $D_{(e)}$ admits a
generalization into the divisor
\vskip3truemm
$$D_{(e)}'\,:\qquad D_0\quad \cou{}{-n}\nlin
\cshiftup{B'_2}{-1}\lin\!\!\lin
\cou{\!\!\!\!\!\!\!\!\!\!\!\! D_1}{-2}\nlin
\cshiftup{B'_1}{-1}\lin \cou{D_2}{-2}
\lin\ldots\lin\cou{D_{n-1}}{-2}\quad.$$

\noindent It is interesting to note that $D'_{(e)}$ does
no longer correspond to
a $\C^*$-surface (see Proposition \ref{eboundary.9} below). \eexas

\section{Extended divisors of $\C^*$-surfaces}

In this section we examine as to when a Gizatullin
$\C^*$-surface has a distinguished or rigid extended divisor. Our
main criterion is Theorem \ref{nonspec} below. We review
first some basic facts about $\C^*$-surfaces.

\subsection{DPD presentation for $\C^*$-surfaces}
\bsit\label{eph} A normal affine surface $V=\Spec A$ endowed with
an effective $\C^*$-action is called a {\it $\C^*$-surface}. Such
a surface can be elliptic, parabolic or hyperbolic. On a non-toric
elliptic or parabolic $\C^*$-surface $V$ the $\C^*$-action is
unique up to conjugation in the automorphism group $\Aut (V)$ and
inversion in $\C^*$ \cite[Corollary 4.3]{FlZa3}. Moreover by
Corollaries 3.23, 4.4 and Theorem 4.5 in \cite{FlZa2} for any
non-toric Gizatullin $\C^*$-surface $V$, the $\C^*$-action on $V$
is hyperbolic. Therefore to deduce Theorem \ref{MT} it is enough
to restrict to hyperbolic $\C^*$-surfaces. \esit

\bsit\label{dpdpr} A simple and convenient description for elliptic
and parabolic $\C^*$-surfaces in terms of the  associated gradings
on the coordinate rings was elaborated by Dolgachev, Pinkham and
Demazure. It was extended to the hyperbolic case in \cite{FlZa1},
where this construction was called a {\it DPD presentation}.

Namely, any hyperbolic $\C^*$-surface $V$
can be presented as
$$
V=\Spec A,\qquad\mbox{where}\quad
A=A_0[D_+,D_-]=A_0[D_+]\oplus_{A_0} A_0[D_-]\,
$$ for a pair  of $\Q$-divisors $(D_+,D_-)$
on a smooth affine curve $C=\Spec A_0$ satisfying the condition
$D_++D_-\le 0$. Here
$$
A_0[D_\pm]=\bigoplus_{k\ge 0} H^0(C,
\cO_C(\lfloor kD_\pm\rfloor))u^{\pm k}
\subseteq \Frac(A_0)[u,u^{-1}]\,,
$$
where $\lfloor D\rfloor$ stands for the integral part of a
divisor $D$ and $u$ is an independent variable. A posteriori,
$u\in\Frac (A_0)\otimes_{A_0}A$ and $\deg\, (u)=1$. One can change
$u$ by multiplying it by a function $\varphi\in\Frac (A_0)$; then
$D_\pm$ will be replaced by $D'_\pm=D_\pm\pm \div\,\varphi$.

We say that two pairs $(D_+,D_-)$ and $(D'_+,D'_-)$ are {\em
equivalent} if $D'_\pm=D_\pm\pm \div\,\varphi$ for a rational
function $\varphi$ on $C$. By Theorem 4.3(b) in \cite{FlZa1} two
hyperbolic $\C^*$-surfaces $V=\Spec A_0[D_+,D_-]$ and $V'=\Spec
A_0[D'_+,D'_-]$ are equivariantly isomorphic over $C=\Spec A_0$ if
and only if the pairs $(D_+,D_-)$ and $(D'_+,D'_-)$ are
equivalent.
\esit

\bsit\label{fxpt} The embedding $A_0\hookrightarrow A_0[D_+,D_-]$
induces an orbit map $\pi:V\to C$. The fixed points on a
hyperbolic $\C^*$-surface $V$ are all isolated, attractive in one
and repelling in the other direction.
The numerical characters of these singular points are precized
in the next result.
\esit

\blem\label{sing} (\cite[Theorem 4.15]{FlZa1}) For a point
$p\in\A^1$ we let \be\label{numerics}
D_+(p)=-\frac{e^+}{m^+}\mbox{ and } D_-(p)=\frac{e^-}{m^-}\mbox{
with } \gcd(e^\pm,m^\pm)=1 \mbox{ and } \pm m^\pm>0\,. \ee
Then the following hold.
\bnum[(a)]
\item If $D_+(p)+D_-(p)=0$ then $\pi^{-1}(p)\cong \C^*$ is a fiber
of multiplicity $m:=m^+=-m^-$ which contains no singular point of
$V$.

\item If $D_+(p)+D_-(p)<0$ then the fiber $\pi^{-1}(p)$ in $V$
consists of two orbit closures $O^\pm\cong \A^1$ of multiplicity
$\pm m^\pm$ in the fiber $\pi^{-1}(p)$ meeting in a unique point $p'$.
Moreover $V$ has a cyclic quotient singularity of type $(\Delta, e)$
at $p'$, where
\be\label{dete} \Delta=\Delta(p')=-\left|\ba{ll}e^+& e^-\\
m^+& m^-\ea\right| =m^+m^-(D_+(p)+D_-(p))>0\,,\ee and $e$ with $0\le
e< \Delta$ is defined by
$$e=e(p')\equiv \left|\ba{ll}a& e^-\\
b& m^-\ea\right| \mod \Delta\quad\mbox{if}\quad
\left|\ba{ll}a& e^+\\
b& m^+\ea\right|=1\,.$$ \enum
\elem

For instance, if $D_\pm(p)$ are both integral
and $k=-(D_+(p)+D_-(p))>0$ then $V$
has an $A_{k-1}$-singularity at $p'$.
We also need the following observation, see \cite[Theorem 4.5]{FlZa2} and
\cite[Lemma 4.2(b)]{FKZ2}.

\blem\label{giz}
For a $\C^*$-surface $V=\Spec A_0[D_+,D_-]$ the following hold.
\bnum[(a)]
\item $V$ is a Gizatullin $\C^*$-surface
if and only if $A_0\cong\C[t]$ and
$\supp\{D_\pm\}\subseteq\{p_\pm\}$
for some points $p_\pm\in\A^1=\Spec \C[t]$.
\item $V$ is toric if and only if $A_0\cong \C[t]$ and up to equivalence
$(D_+,D_-)$ is the divisor $(\frac{-e^+}{m^+}[p_0],
\frac{e^-}{m^-}[p_0])$, for some point $p_0\in\A^1=\Spec \C[t]$.
\enum \elem

\subsection{Completions of $\C^*$-surfaces}
We let $V=\Spec A_0[D_+,D_-]$ be a normal affine $\C^*$- surface. We
review here some facts on equivariant completions of $V$; for proofs
we refer the reader to \cite{FKZ2}.

\blem \label{barcomp} $V$ admits an equivariant  normal completion
$(\bX,\bD)$ with the following properties. \bnum[1.] \item  {\em
(Cf.\ \cite[Proposition 3.8 and Remark 3.9(4)]{FKZ2})} The orbit
map $V\to C=\Spec A_0$ extends to a $\PP^1$-fibration $\pi: \bX\to
\bC $, where $\bC$ is the smooth completion of $C$. \item  $\bD$
has exactly two horizontal components $\bC_\pm$, which are
sections of $\pi$, where $\bC_+$ is repelling and $\bC_-$ is
attractive. \item {\em (Cf.\ \cite[3.10 and Proposition
4.18]{FlZa2})} For $D_+(p)+D_-(p)=0$ the fiber
$\bO_p=\pi^{-1}(p)\cong \PP^1$ has multiplicity $m^+=-m^-$, where
$m^\pm$ are as in (\ref{numerics}). \item  {\em (Cf.\ \cite[3.10
and Proposition 3.13(d) ]{FKZ2} and \cite [Proposition
4.18]{FlZa2})} If $D_+(p)+D_-(p)<0$ then\footnote{with the
notation as in Lemma \ref{sing}.} the fiber $\pi^{-1} (p)$
consists of two orbit closures $\bO_p^\pm\cong \PP^1$ of
multiplicity $\pm m^\pm$ meeting in a unique point $p'$ (cf.\
Lemma \ref{sing}(b)). Moreover $\bO_p^\pm$ have self-intersection
indices $\frac{m^\mp}{\Delta m^\pm}$, respectively. \enum \elem

In general, $\bD$ can contain singular points of $\bX$.
Let $\rho:\tX\to\bX$ denote the minimal resolution of singularities
of $\bX$ and $\tD:=\rho^{-1}(\bD)$.
The $\C^*$-action on $\bX$ then lifts to $\tX$.

\blem\label{resolution} {\rm
(See \cite[Proposition 3.16]{FKZ2})}
Let $\tilde\pi:\tX\to \bC$ be the induced
$\PP^1$-fibration and let
$\tC_\pm$ be the proper transforms of $\bC_\pm$ on $\tX$. Then the
following hold.
\begin{enumerate}[(a)]
\item $(\tX, D)$ is an SNC completion of the
minimal resolution of $V'$ of $V$.  Moreover, $\tC_\pm^2=\deg \lfloor
D_\pm\rfloor $.

\item If $D_+(p)+D_-(p)<0$ then the fiber
$\tilde\pi^{-1}(p)$ together with $\tC_\pm$ has dual graph
\vspace*{0.3truecm} \be\label{dovos} \co{\tC_+}
\vlin{15}\boxo{\{D_+(p)\}} \vlin{15}\co{\tO^+_p}
\vlin{15}\boxo{\left(e/ \Delta\right)^*} \vlin{15}\co{\tO^-_p}
\vlin{15}\boxo{\{D_-(p)\}^*} \vlin{15}\co{\tC_-}\qquad, \ee where
$\tO^\pm_p$ with $(\tO^\pm_p)^2=\lfloor\frac{m^\mp}{\Delta
m^\pm}\rfloor$ are the proper transforms of $\bO_p^\pm$,
respectively, and at least one of them is a $(-1)$-curve
\footnote{See \ref{not1} for the notation $(e/\Delta)^*$.}.

\item If $D_+(p)+D_-(p)=0$ then the fiber $\tilde\pi^{-1}(p)$
together with $\tC_\pm$ has dual graph \vspace*{0.3truecm}
\be\label{+=-} \co{\tC_+} \vlin{15}\boxo{\{D_+(p)\}}
\vlin{15}\co{\tO_p} \vlin{15}\boxo{\{D_-(p)\}^*}
\vlin{15}\co{\tC_-}\qquad, \ee where the proper transform $\tO_p$
of $\bO_p$ is a $(-1)$-curve.
\end{enumerate}
\elem

\brem\label{fibreak} 1. If $D_\pm (p)\in\Z$ and
$-(D_+(p)+D_-(p))=\Delta>0$, then by Lemma \ref{resolution}(b) $V$
has an $A_{\Delta-1}$-singularity at $p'$ and the graph
(\ref{dovos}) is
\be\label{dovo} \,\qquad\co{\tC_+} \vlin{15}\cou{\tO^+_p}{-1}
\vlin{15}\boxo{A_{\Delta-1}} \vlin{15}\cou{\tO^-_p}{-1}
\vlin{15}\co{\tC_-}\qquad. \ee \erem

\subsection{Extended divisors of Gizatullin $\C^*$-surfaces}
\bsit\label{sit33} In this section we let $V=\Spec \C[t][D_+,D_-]$
denote a Gizatullin $\C^*$-surface. By  Lemma \ref{giz}(a) $\supp
\{D_+\}\subseteq \{p_+\}$, $\supp \{D_-\}\subseteq \{p_-\}$,
 and the orbit map $\tilde\pi:\tX\to\bC=\PP^1$
is defined by the linear system $|F_\infty|$, where $F_\infty$
denotes the fiber $\tilde\pi^{-1}(\infty)$ over the point
$\{\infty\}=\PP^1\setminus\A^1$. Furthermore by Lemma
\ref{resolution} the boundary zigzag $\tD$ of $V$ in $\tX$ has
dual graph \be \label{dovos1}
\tD:\qquad\qquad\boxo{\{D_+(p_+)\}^*} \vlin{15}\cou{\tC_+}{}
\vlin{15}\cou{F_\infty}{0} \vlin{15}\cou{\tC_-}{}
\vlin{15}\boxo{\{D_-(p_-)\}} \qquad. \ee A standard equivariant
completion $\tV$ of the resolution $V'$ of $V$ can be obtained
from $\tX$ by moving the zero weight in (\ref{dovos1}) to the left
via elementary transformations \cite{FKZ1}. More precisely, to
obtain the standard zigzag $D=C_0+\ldots+C_n$ from (\ref{dovos1})
one has to perform first a sequence of elementary transformations
at $F_\infty$ until $\tC_+$ becomes a $0$-curve. At this step the
self-intersection index of the image $C_s$ of $\tC_-$ becomes
equal to $ w_s=\deg\, (\lfloor D_+\rfloor+\lfloor D_-\rfloor)$. By
moving the two resulting neighboring zeros to the left  via
a sequence of elementary transformations (which contracts in
general the curve $\tC_+$ and does not affect $C_s^2=w_s$) one
gets a completion $(\tV,D)$ of $V'$ by the zigzag \be\label{bo}
D:\qquad\qquad \cou{C_{0}}{0} \lin\cou{C_{1}}{0} \vlin{15}
\boxo{\{D_+(p_+)\}^*} \vlin{15}\cou{C_s}{w_s}
\vlin{15}\boxo{\{D_-(p_-)\}} \qquad.\ee This zigzag is standard as
soon as $w_s\le -2$. Indeed, all curves in the boxes labelled
$\{D_+(p_+)\}^*$ and $\{D_-(p_-)\}$ have weight $\le -2$. The
elementary transformations as above result in a birational
morphism $\tX\dashrightarrow\tV$, which is  the identity on $V'$.

By abuse of notation we keep the same symbols $\tO_p$,
$\tO^-_{p_\pm}$ in both completions $\tX$ and $\tV$, cf.
(\ref{dovos}). Note that the self-intersection indices
$(\tO^-_{p_\pm})^2$, $\tO_p^2$ are the same in $\tX$ and in
$\tV$. \esit

To describe the resulting extended graph it is convenient to
introduce {\em admissible} feather collections
$\{\fF_\rho\}_{\rho\ge 1}$; see \cite{FKZ2}. By this we mean that
all but at most one feather $\fF_\rho$ are $A_k$-feathers. Further,
a curve on a $\C^*$-surface is called {\em parabolic} if it is
pointwise fixed. For the next result, we refer the reader to
Proposition 5.8 in \cite{FKZ2} and its proof.

\bprop\label{eboundary.9} With the notations as above, the
resolution of singularities $V'\to V$ of a normal affine Gizatullin
$\C^*$-surface $V$ admits an equivariant SNC completion $(\tV, D)$
with extended graph

\vskip0.3truecm

\be\label{ezigzag1} D_{\rm ext}: \qquad\quad\cou{C_0}{0}\lin
\cou{C_1}{0}\vlin{18} \boxou{}{\{D_+(p_+)\}^*} \vlin{18}
\cou{C_s\quad\;\;}\, \cou{}{w_s}\nlin\xbshiftup{}{\qquad
\{\fF_\rho\}_{\rho \ge 1}}\vlin{20} \boxou{}{\{D_-(p_-)\}}
\nlin\boxshiftup{}{ \fF_0} \quad\quad\ee and with boundary zigzag
$D$ represented by the bottom line in (\ref{ezigzag1}). Here
$w_s=\deg (\lfloor D_+ \rfloor+\lfloor D_-\rfloor)$, $\fF_0$ is a
single feather (possibly empty), $\{\fF_\rho\}_{\rho\ge 1}$ is an
admissible feather collection and $C_s=\tC_-$ is an attractive
parabolic component. Moreover the following hold:
\begin{enumerate} [(a)] \item
The feather collection $\{\fF_\rho\}_{\rho\ge 1}$ is empty if and
only if $V$ is a toric surface. If $V$ is
non-toric\footnote{However, see Remark \ref{vosst}(4) below.} then
$w_s\le -2$ and consequently $(\tV, D)$ is a standard completion
of $V'$.
\item If $p_+\ne p_-$ then $(\tO_{p_-}^-)^2=-1$ and the
feathers \be\label{sewe1} \fF_0:\qquad \cou{\tO_{p_-}^-}{}
\vlin{20}\boxo{(e/\Delta)(p_-)} \qquad\qquad
\mbox{and}\qquad\qquad \fF_1:\qquad
 \cou{\tO_{p_+}^-}{}
\vlin{20}\boxo{(e/\Delta)(p_+)}\quad\quad \ee are contained in
the fibers over $p_-$ and $p_+$, respectively, as described in
(\ref{dovos}). \item If $p_+=p_-=:p$ then the $\fF_\rho$ are $A_
{k_\rho}$-feathers $\forall \rho\ge 1$.
 The feather $\fF_0$
 is empty
 if and only if $D_+(p)+D_-(p)= 0$. Otherwise it is as in
 (\ref{sewe1})
 with $p_-=p$ and   $(\tO_{p_-}^-)^2=\lfloor\frac{m^+}{\Delta
m^-}\rfloor$.
\end{enumerate}
\eprop

\bproof By virtue of Lemma 2.20 in \cite{FKZ2}, $V$ is toric
if and only if the extended divisor $D_\ext$ is linear. This
yields the first assertion in (a). Thus by Proposition 5.8 in
\cite{FKZ2}, only the second assertion in (a) and the first one in
(b) need to be proved.

Assuming that $w_s=0$ it is easily seen that $(D_+,D_-) \sim (0,0)$,
$D=[[0,0,0]]$ and $\fF_\rho=\emptyset\,\,\forall \rho\ge 0$. But
then $V\cong\A^1\times\C^*$, which contradicts our assumption that
$V$ is a Gizatullin surface. Thus $w_s\le-1$. If $w_s=-1$ then
necessarily $p_+=p_-$ and $\lfloor D_+(q)\rfloor +\lfloor
D_-(q)\rfloor =0$ for any point $q$ different from $p:=p_+=p_-$.
Since $D_+(q)+D_-(q)\le 0$ it follows that $D_+(q)=-D_-(q)$ are
integral for $q\ne p$. Passing to an equivalent pair of divisors we
may suppose that $D_+$ and $D_-$ are both supported at $p$. Hence
$V$ is toric by Lemma 4.2(b) in \cite{FKZ2}\footnote{Cf. Claim
($\alpha$) in the proof of Proposition 5.8 in \cite{FKZ2}.}.

Finally, the equality $(\tO_{p_-}^-)^2=-1$ in (b) follows from
Lemma \ref{resolution}(b). Indeed, as $D_+(p_-)\in\Z$ we have
$m^+=m^+(p_-)=1$ and so $(\tO_{p_-}^-)^2=\lfloor\frac{m^+}{\Delta
m^-}\rfloor=-1$.\eproof

\brems\label{vosst} 1. One can also move the zeros in
(\ref{dovos1}) to the right. In the case where $w_s\le -2$ this
yields a second standard completion with the boundary zigzag
reversed. However, in this completion $(\tV^\vee,D^\vee)$ the
parabolic component is repelling, and it becomes attractive when
replacing the given $\C^*$-action by the inverse one via the
automorphism $t\mapsto t^{-1}$ of $\C^*$. The extended dual graph
$D_\ext$ in (\ref{ezigzag1}) is uniquely determined by the
requirement that it corresponds to an equivariant standard
completion of $V'$ with attractive parabolic component.

2. If $V$ is smooth then every feather in (\ref{ezigzag1}) consists
of a single irreducible curve, see \ref{smooth}. A more detailed
description can be found in \cite[Corollary 5.10]{FKZ2}. If for
instance $p_+\ne p_-$ or one of the fractional parts $\{D_+\}$,
$\{D_-\}$ vanishes then, up to passing to an equivalent pair of
$\Q$-divisors,
$$\left(D_+,D_-\right)=
\left(-\frac{1}{k}[p_+]\,\,, -\frac{1}{l}[p_-]-D_0\right)
\qquad\mbox{with}\quad k,l\ge 1 \,,$$ where $D_0=\sum_{\rho=2}^t
[p_\rho]$ is a reduced integral divisor on $C\cong \A^1$ so that all
points $p_\rho$ are pairwise distinct and different from $p_\pm$.
Thus in (\ref {ezigzag1}) the boxes adorned $ \{D_+(p_+)\}^*$ and
$\{D_-(p_-)\}$ are just $A_{k-1}$- and $A_{l-1}$-boxes, which
represent chains of $(-2)$-curves $[[(-2)_{k-1}]]$ and $
[[(-2)_{l-1}]]$, respectively.

3. Contracting the exceptional curves in $\tV$ corresponding to the
singularities in the affine part $V$  we obtain a standard completion
$(\bV,D)$ of $V$.

4. For a toric Gizatullin surface it may happen that $w_s= -1$,
take e.g.\ $V=\Spec \C[t][D_+,D_-]$ with
$(D_+,D_-)=\left(-\frac{1}{2}[0], \,\frac{1}{3}[0]\right)$. The
boundary zigzag as in Proposition \ref{eboundary.9} is now
$[[0,0,-2,-1,-3]]$, which has standard form $[[0,0]]$. Thus
$V\cong \A^2$.\erems

Let us compute more generally the standard boundary zigzag of
arbitrary affine toric surface $V=V_{d,e}=\A^2/\Z_d$ (see
\ref{toricsit}), where $0\le e<d$ and $\gcd (e,d)=1$.

\blem\la{toric}
The toric surface $V_{d,e}$ admits a standard completion with
boundary zigzag
\be\label{to*} D:\qquad\quad\co{0}\lin\co{0}\llin \boxo{
\frac{d-e}{d}}\;\;.
\ee
Moreover, the reverse zigzag $D^\vee$ is given by
$\quad\co{0}\lin\co{0}\llin \boxo{\frac{d-e'}{d}}\quad $,
where $e'$ is the unique number with $0\le e'< d$ and
$ee'\equiv 1\mod d$. In particular, the standard boundary
of a toric surface is symmetric if and only if $e^2\equiv 1 \mod d$.
\elem

\bproof
Using Lemmas \ref{sing}
and \ref{giz}(b), $V_{d,e}\cong \Spec \C[t][D_+,D_-]$ with $D_+=0$
and $D_-=\frac{d}{e-d}[0]$. According to Proposition \ref{eboundary.9} the
standard boundary has dual graph
$$
\co{0}\lin\co{0}\lin\co{\lfloor
\frac{d}{e-d}\rfloor} \llin\boxo{\{ \frac{d}{e-d}\}}\;\;.
$$
A simple computation gives
$$
\co{\lfloor \frac{d}{e-d}\rfloor} \llin\boxo{\{ \frac{d}{e-d}\}}
\qquad =\qquad \boxo{ \frac{d-e}{d}}\quad.
$$
Finally, the form of $D^\vee$ follows from \ref{not1}.
\eproof

\brem\label{to}
1. The form of $D^\vee$ reflects
the well known
fact that $V_{d,e}\cong V_{d',e'}$ if and only if $d=d'$ and either
$e=e'$ or $ee'\equiv 1 \mod d$, see e.g. \cite[Remark 2.5]{FlZa1}.

2. Due to the lemma, the toric surface
$V_{d,e}$ is uniquely determined by its standard
boundary zigzag.
\erem

For later use we give a criterion as to when a $\C^*$-action is
equivalent to its inverse.

\blem\label{symmetric0} For a $\C^*$-surface $V=\Spec A_0[D_+,D_-]$
over $C=\Spec \,A_0$, the  associated  hyperbolic $\C^*$-action
$\Lambda$ on $V$ and its inverse action $\Lambda^{-1}$ are conjugate
in the automorphism group $\Aut (V)$ if and only if there exists an
automorphism $\psi\in\Aut (C)$ such that \bnum[(i)]\item
$\psi^*(D_++D_-)=D_++D_-$ and\item $\psi^*(D_-)-D_+$ is a principal
divisor. \enum \elem

\bproof Inverting the $\C^*$-action results in interchanging the
components $A_0[D_+]$ and $A_0[D_-]$ of the graded algebra
$A_0[D_+,D_-]$ or, equivalently, in interchanging the divisors $D_+$
and $D_-$ (see Section 3.1). Thus the inverse action $\Lambda^{-1}$
corresponds to the $\C^*$-surface $V^\vee=\Spec A_0[D_-,D_+]$ over
$C$. By Theorem 4.3(b) in \cite{FlZa1}, the actions $\Lambda$ and
$\Lambda^{-1}$ are conjugate in the group $\Aut\, V$ if and only if
the $\C^*$-surfaces $(V,\Lambda)$ and $(V^\vee,\Lambda^{-1})$ are
equivariantly isomorphic, if and only if there is an automorphism,
say, $\psi$ of $C$ such that the pairs $(D_+,D_-)$ and
$(\psi^*(D_-), \psi^*(D_+))$ are equivalent i.e.,
$$D_++D_0=\psi^*(D_-)\qquad\mbox{and}\qquad
D_--D_0=\psi^*(D_+)\,$$ for some principal divisor $D_0$ on $C$. The
first of these equalities yields (ii), and taking their sum gives
(i). \eproof

\brems\label{srt}
1. Suppose that $V=\Spec \C[t][D_+,D_-]$ is a $\C^*$-surface over
$\A^1=\Spec \,\C[t]$. Then condition (ii) in Lemma \ref{symmetric0}
is equivalent to
$\psi^*(\{D_+\})=\{D_-\}$. In particular, if the divisor $D_+-D_-$
is integral then (i) and
(ii) are automatically satisfied with $\psi=\id$.

2. We have seen in Remark \ref{vosst}(1) that changing the
$\C^*$-action of a Gizatullin $\C^*$-surface $V=\Spec
\C[t][D_+,D_-]$ by the automorphism $t\mapsto t^{-1}$ of $\C^*$
amounts to reversing the standard zigzag. So {\em if the
$\C^*$-action on $V$ is conjugate to its inverse then the standard
zigzag $D$ of $V$ is symmetric.}

3. Note however that for a Gizatullin $\C^*$-surface with a
symmetric standard boundary zigzag the $\C^*$-action is not
conjugate to its inverse, in general. A simple example is given by the toric surface $V=\Spec\C[t][D_+,D_-]\cong \A^2$ with
$(D_+,D_-)=\left(-\frac{1}{2}[0],\, \frac{1}{3}[0]\right)$, see Remark \ref{vosst}(4). This pair does not satisfy condition (ii) of Lemma
\ref{symmetric0} although its standard boundary
zigzag is equal to $[[0,0]]$ and so is symmetric.
\erems

\subsection{A rigidity criterion}
In Theorem \ref{nonspec} below we show that under the assumptions
($\alpha_+$) and ($\beta$) of Theorem \ref{MT}, the standard
divisor (\ref{extended1}) is distinguished and rigid. Moreover, if
($\alpha_*$) holds then this divisor is rigid after
possibly interchanging $D_+$ and $D_-$.

\bsit\label{greek} We begin by recalling the assumptions
($\alpha_+$), ($\alpha_*$) and ($\beta$) of Theorem
\ref{MT}.

\begin{enumerate}
\item[($\alpha_+$)] $\supp \{D_+\}\cup \supp \{D_-\}$ is empty or
consists of one point, say, $p$ satisfying either
$D_+(p)+D_-(p)=0\,$  or
\be\label{alphaprim}
D_+(p)+D_-(p)\le -\max\left(\frac{1}{{m^+}^2},\,
\frac{1}{{m^-}^2}\right)\,,\ee where $\pm m^\pm$ is the minimal
positive integer such that $m^\pm D_\pm(p)\in\Z$.

\item[($\alpha_*$)] $\supp \{D_+\}\cup \supp \{D_-\}$ is empty or
consists of one point  $p$, where
$$
D_+(p) + D_-(p) \le -1
\quad\mbox{or}\quad \{D_+(p)\} \ne 0\ne \{ D_-(p)\}.
$$

\item[($\beta$)] $\supp \{D_+\}=\{p_+\}$ and $\supp
\{D_-\}=\{p_-\}$ for two different points $p_+,p_-$, where
\be\label{betaprim} D_+(p_+) + D_-(p_+)\le-1\quad\mbox{and}\quad
D_+(p_-) + D_-(p_-)\le-1 \,. \ee
\enum
\esit

\blem \label{-1} For a point  $p\in\A^1$  with $(D_++D_-)(p)<0$ the
following hold.
\bnum[(a)]\item
$\tO^\pm_p$ in (\ref{dovos}) is a $(-1)$-curve if and only if $(D_++D_-)(p)\le -1/(m^\pm)^2$. In particular,
both $\tO^+_p$ and $\tO^-_p$ in
(\ref{dovos}) are $(-1)$-curves\footnote{Anyway, at least one of
these is a $(-1)$-curve, see Lemma \ref{resolution}(b).} if and only
if (\ref{alphaprim}) is fulfilled.
\item If $\min\left(\{D_+(p)\},\, \{D_-(p)\}\right)=0$ then
(\ref{alphaprim}) is equivalent to
\be\label{-2e} D_+(p)+D_-(p)\le -1\,.\ee \enum
\elem

\bproof We let as before $D_\pm(p)=e^\pm/m^\pm$ with $\gcd(e^\pm,
m^\pm)=1$, $m^+,-m^-\ge 1$ and
$$\Delta=\Delta(p)=m^+m^-(D_+(p)+D_-(p))\ge 1\,.$$
(a) follows from the equalities
$(\tO^\pm_p)^2=\lfloor\frac{m^\mp}{\Delta m^\pm}\rfloor$, see
Lemma \ref{resolution}(b). Indeed,
$$\left\lfloor\frac{m^\mp}{\Delta m^\pm}\right\rfloor=-1\;
\Longleftrightarrow\;\frac{m^\mp}{\Delta m^\pm}\ge-1\;\Longleftrightarrow\;\frac{m^\mp}{m^\pm}\ge-\Delta\;
\Longleftrightarrow\;\frac{-1}{(m^\pm)^2}\ge(D_++D_-)(p).$$

To show (b), after interchanging $D_+$ and $D_-$, if necessary, and
passing to an equivalent pair of divisors, which does not affect our
assumptions, we may suppose that $D_+(p)=0$. Thus $m^-\le -1$ and
$m^+=1$ and so
$$\max\left(\frac{1}{{m^+}^2},\,\frac{1}{{m^-}^2}\right)
=\max\left(1,\,\frac{1}{{m^-}^2}\right)=1\,.$$ Now (b) follows.
\eproof

For a Gizatullin $\C^* $-surface $V=\Spec \C[t][D_+,D_-]$ we let
$V^\vee=\Spec \C[t][D_-,D_+]$, and we denote by
$D_\ext,\,D_\ext^\vee$ the corresponding extended divisors of the
equivariant standard completions as in Proposition
\ref{eboundary.9}.

\blem\label{bridge-1}
\bnum[(a)]\item All bridges $B_\rho $ of the
feathers $\fF_\rho $, $\rho\ge 0$, are $(-1)$-curves in both
divisors $D_\ext$ and $D_\ext^\vee$ if and only if ($\alpha_+$) or
($\beta$) holds.
\item If $(\alpha_*)$ is
fulfilled then all bridges $B_\rho $ are  $(-1)$-curves in at
least one of these divisors.
\enum \elem

\bproof
Assume first that $p_+=p_-=p$. By
Proposition \ref{eboundary.9}(c) the $\fF_\rho$  are
$A_{k_\rho}$-feathers $\forall \rho\ge 1$. Hence the corresponding
bridges are $(-1)$-curves. If $D_+(p)+D_-(p)= 0$ then again by
Proposition \ref {eboundary.9}(c), $\fF_0=\emptyset$ and we are
done. If $D_+(p)+D_-(p)< 0$ then by Lemma \ref{-1} the remaining
bridges $\tO^\pm_p$ of the feather $\fF_0$ in both $D_\ext$ and
$D_\ext^\vee$ are $(-1)$-curves if and only if (\ref{alphaprim})
holds, as claimed in (a). Anyhow, according to Lemma
\ref{resolution}(b) at least one of $\tO_p^\pm$ is a $(-1)$-curve,
hence (b) follows as well in this case.

Suppose further that $p_+\neq p_-$.
By Proposition \ref{eboundary.9}(b) the bridge $\tO_{p_-}^-$ of
the feather $\fF_0$ in $D_\ext$ is a $(-1)$-curve and, symmetrically, the bridge $\tO_{p_+}^+$ of the feather $\fF_0$ in $D_\ext^\vee$ is a
$(-1)$-curve.
Thus by Lemma \ref{-1} the bridge $\tO_{p_+}^-$  of the
feather $\fF_1$ in $D_\ext$\footnote{See (\ref{sewe1}).} is a $(-1)$-curve if
and only if the first inequality in (\ref{betaprim})
is fulfilled.  Similarly
the bridge $\tO_{p_-}^+$  of the
feather $\fF_1$ in $D_\ext^\vee$ is a $(-1)$-curve if
and only if the second inequality in (\ref{betaprim})
is satisfied. The other bridges are as well $(-1)$-curves due to
the fact that the feather collection $\{\fF_\rho\}$ is admissible,
see Proposition \ref{eboundary.9}. This implies (a)
in this case.
\eproof

\brem\label{interch} Switching $D_+$ and $D_-$ amounts to
interchanging $D_\ext$ and $D_\ext^\vee$. So replacing the given
$\C^*$-action by its inverse one can achieve, if necessary, that
the conclusion of Lemma \ref{bridge-1}(b) holds for the model with
an attractive parabolic component.\erem

\bdefi\label{tail1} Suppose that $\supp \{D_+\}\subseteq \{p_+\}$
and $\supp \{D_-\}\subseteq \{p_-\}$ with (not necessarily
distinct) points $p_\pm$. By the {\em tail} of the extended
divisor (\ref {ezigzag1}) we mean the subgraph
\be\label{tail}
L=L_{s+1}=\qquad \quad \boxo{\{D_-(p_-)\}} \vlin{17}\boxo{ \fF_0}
\qquad\mbox{=} \qquad \co{C_{s+1}} \lin \ldots \lin
\co{C_n}\llin\cou{\tO^-_{p_-}}{}\vlin{15}
\boxo{(e/\Delta)(p_-)}\qquad,
\ee
cf. (\ref{ezigzag1}), (\ref{sewe1}),
and by a {\em subtail} a subgraph of $L$ of
the form
\be\label{subtail} L_t=\qquad \quad \co{C_{t}} \lin
\ldots \lin \co{C_n}\llin\cou{\tO^-_{p_-}}{}\vlin{15}
\boxo{(e/\Delta)(p_-)} \ee with $s+1\le t\le n$.
\edefi

\blem \label{lemtail0}
If $D_+(p_-)+D_-(p_-)\ne 0$ then
the tail $L$ is contractible if and only if $\{D_+(p_-)\}=0$.
In particular, if $p_+\neq p_-$ then $L$ is contractible\footnote{Cf. \cite[Proposition
5.8]{FKZ2}.}.
\elem

\bproof Suppose first that $L$ is contractible. By Lemma
\ref{resolution}(b) the fiber $\tilde\pi^{-1}(p)$ with $p:=p_-$
has dual graph \vspace*{0.3truecm} \be\label{onlyfiber}
\boxo{\{D_-(p)\}} \llin\co{\tO^-_p} \llin\boxo{e/ \Delta}
\llin\co{\tO^+_p} \llin\boxo{\{D_+(p)\}^*} \qquad=\qquad \boxo{L}
\llin\co{\tO^+_p} \llin\boxo{\{D_+(p)\}^*} \qquad , \ee where we
use the notations of {\em loc.cit.}. If $L$ is contractible then
contracting it in the fibre (\ref{onlyfiber}) leads to the
divisor\quad $\co{A}\llin\boxo{\{D_+(p)\}^*} $\qquad, where
$A$ denotes the image  of $\tO_p^+$ and all the weights in the box
adorned $\{D_+(p)\}^*$ are $\le -2$. This divisor has to be
contractible to a smooth fibre $[[0]]$, which is only possible if
the box is empty.

Conversely, if  $\{D_+(p)\}=0$ then by \ref{barcomp}(4) and Lemma
\ref{resolution}(b)
 $\tO_{p}^+$ has multiplicity 1 in the fibre
(\ref{onlyfiber}), hence the rest of it, which is $L$, can be
contracted to a smooth point. \eproof

\blem \label{lemtail}
 \bnum[(a)]
\item If $p_+\ne p_-$ and  $(D_++D_-)(p_-)\le -1$ then none of the
subtails $L_t$ with $t\ge s+2$ is contractible. The same holds if
$(\alpha_+)$ is satisfied. \item If $p_+\ne p_-$ and
$0>(D_++D_-)(p_-)>-1$ then the subtail $L_{s+2}$ is contractible.
\item If $p_+= p_-=:p$, $(D_++D_-)(p)\ne 0$ and (\ref{alphaprim})
is not satisfied then either $\tO^-_p$ is not a $(-1)$-curve or
the subtail $L_{s+2}$ is contractible. \enum \elem

\bproof If in (a) $(D_++D_-)(p_-)=0$ then $\fF_0=\emptyset$ and so
every non-empty subtail of $L$ is minimal and hence
non-contractible. Otherwise $\fF_0\ne\emptyset$, and under the
assumptions of (a) Lemma \ref{-1}(a) implies $(\tO_{p_-}^+)^2=-1$.
If a proper subtail $L_t$ of $L$ were contractible then, while
contracting the fiber (\ref{onlyfiber}) with $p=p_-$ to $[[0]]$,
at least one component neighboring $\tO_{p_-}^+$ would be
contracted. Hence the image of $\tO_{p_-}^+$ would have
self-intersection $\ge 0$ and so it must be the full fiber. This
contradicts the assumption that $t\ge s+2$ and so (a) holds.

(b) In this case $(\tO_{p_-}^+)^2\le -2$, see Lemma \ref{-1}(a). If
(b) does not hold then contracting $L$, $C_{s+1}$ must be
contracted before the subtail $L_{s+2}$ is contracted. It follows
that there is a proper contractible subchain, say, $P$ of $L$
which contains the piece $[C_{s+1}, \ldots, C_n, \tO_{p_-}^-]$.
Contracting $P$ in the full fiber (\ref {onlyfiber}) leads to a
linear chain \be\label{oo} \co{-1} \lin\co{E_1}\lin \ldots\lin
\co{E_s}\lin \co{\tO^+_{p_-}} \qquad, \ee where $[E_1,\ldots,
E_s]$ is a subchain of the box labelled by $e/\Delta$. However,
since all curves in $[E_1,\ldots, E_s, \tO^ +_{p_-}]$ have
self-intersection $\le -2$, (\ref{oo}) cannot be blown down to
$[[0]]$, which gives a contradiction.

(c) By Lemma \ref{-1}(a) one of the curves $\tO^\pm_p$ is not a $(-1)$-curve. Thus, if $\tO^-_p$ is a $(-1)$-curve then $(\tO_{p_-}^+)^2\le -2$. Arguing as in (b) it follows that the subtail $L_{s+2}$ is contractible.
\eproof

\blem\label{disti}
Suppose that ($\alpha_+$) or ($\beta$) holds.
Then the divisors $D_\ext,\,D_\ext^\vee$ are both
distinguished\footnote{See Definition \ref{dist}.}.
\elem

\bproof Since the conditions ($\alpha_+$) and ($\beta$) are
symmetric in $D_+$, $D_-$, it suffices to show that $D_\ext$ is
distinguished. If for some $i$ with $3\le i\le s$ the divisor
$D_\ext^{>i}=D_{(e)}^{>i-2}$ were contractible then after
contracting $D_{(e)}^{>i-2}$ inside $D_{(e)}$ we would obtain as
dual graph \be\label{suchdiv} \cou{C_2}{w_2}\lin\ldots \lin
\cou{C_{i-1}}{w_{i-1}}\lin
\cou{C_i}{-1}\quad,\qquad\mbox{where}\qquad w_j\le -2\quad\forall
j=2,\ldots, i-1\,. \ee However, $D_{(e)}$ can be contracted to
$[[0]]$ while (\ref{suchdiv}) cannot, a contradiction. Thus it is
enough to consider the divisors $D_\ext^{>i}$ with $i\ge s+1$.

If $(\alpha_+)$ or ($\beta$) holds then  by Lemma \ref{lemtail}(a) the
divisors $D_\ext^{>i}$ are not contractible for all $i=s+1,
\ldots, n$. Therefore $D_\ext$ is distinguished.
\eproof

\bthm\label{nonspec}
If $V=\Spec\C[t][D_+,D_-]$ is a Gizatullin $\C^* $-surface then the following hold.
\bnum
\item If ($\alpha_+$) or ($\beta$)  is fulfilled
then both divisors $D_\ext$, $D_\ext^\vee$ are distinguished and
rigid.
\item If ($\alpha_*$) holds
then at least one of the divisors $D_\ext$, $D_\ext^\vee$ is
rigid.
\enum
\ethm

\bproof
(1) Since the conditions ($\alpha_+$) and
($\beta$) are stable under interchanging
$D_+$ and $D_-$, it is enough to consider the extended divisor
$D_\ext$ for the standard completion of $V$. By Lemmas
\ref{bridge-1}(a) and \ref{disti}
$D_\ext$ is distinguished and all its bridges are
$(-1)$-curves. In particular, no feather can jump to the left, see
Proposition \ref{corgen}.

If the feather collection $\{\fF_{sj}\}$ is empty then also no
feather can jump to the right, so $D_\ext$ is rigid. Moreover,
$D_\ext$ is rigid if one of the conditions (i), (ii) of Theorem
\ref{maincrit} is fulfilled.

Suppose further that $\{\fF_{sj}\}\neq\emptyset$ but \ref{maincrit}(i) fails. Then $s<n$ and $D_\ext^{>n}=\emptyset$. In particular
$\fF_0=\emptyset$, and so by Proposition \ref{eboundary.9}(b,c)
$p_+=p_-=:p$ and $D_+(p)+D_-(p)=0$. Since $s<n$ and
$\fF_0=\emptyset$ the tail $L=D_\ext^{\ge s+1}$ is non-empty and
contains only curves of self-intersection $\le -2$. Thus $L$
cannot be contractible and so \ref{maincrit}(ii) holds, whence (1) follows.

(2) In view of (1) we have to consider only the case that
$\{D_+(p)\}\ne 0$, $\{D_-(p)\}\ne 0$ and $D_+(p)+D_-(p)\ne 0$.
By Lemma \ref{bridge-1}(b), after interchanging $D_\pm$
if necessary, the bridge curves of the extended divisor
$D_\ext$ are all $(-1)$-curves. In  particular, no feather can jump to the left. According to Lemma \ref{lemtail0} the tail $L$ is not contractible and so condition (c) in Proposition \ref{corspec} is violated. Thus none of the feathers $\fF_{s\rho}$ can jump to the right and $D_\ext$ is rigid as required.
\eproof

\brem\label{bothnot}
1. It is worthwhile to remark that Theorem  \ref{nonspec}(1)
is sharp. More precisely, let us establish the following.
\begin{enumerate}[(a)]
\item If neither $(\alpha_*)$ nor $(\beta)$ are satisfied then
none of the divisors $D_\ext$,  $D_\ext^\vee$ is rigid. If
$\supp\{D_+\}\cup \supp\{D_+\}$ consists of two distinct points
and $(\beta)$ is violated then at least one of them is not
distinguished.
\item If $p_+=p_-=p$ and $(\alpha_+)$ fails then
none of the divisors $D_\ext$, $D_\ext^\vee$ is at the same time
distinguished and rigid.
\end{enumerate}

\bproof Let us first deduce (a) in the case $p_+\ne p_-$. This
means that $\{D_+(p_+)\}$, $\{D_-(p_-)\}\ne 0$ while one of the
two numbers $(D_++D_-)(p_\pm)$ is $>-1$. By symmetry it suffices
to show that $D_\ext$ is non-rigid.

If $(D_++D_-)(p_+)> -1$ then the bridge $\tO_{p_+}^-$ of the
feather $\fF_1$ for $V$ (cf.\ (\ref{sewe1})) has self-intersection
$\le -2$ (see Lemma \ref{-1}(a)) and so $D_\ext$ is non-rigid. If
$(D_++D_-)(p_-)>-1$ then by Lemmas \ref{lemtail0} and
\ref{lemtail}(b) the tail $L=L_{s+1}$ in (\ref{tail}) and its
subtail $L_{s+2}$ are both contractible. In other words, the
divisors $D_\ext^{> s+1}$ and $D_\ext^{\ge s+1}$ are both
contractible. Thus by Proposition \ref{convspec} with $i=s$ and
$i'=s+1$ any feather $\fF_\rho=\fF_{s,\rho}$, $\rho \ge 1$, can
jump to a feather $\fF'_{s+1,\rho'}$ under a suitable
specialization, and again $D_\ext$ is non-rigid. Moreover it
is non-distinguished. By interchanging $D_+$ and $D_-$, if
necessary, the assumption $(D_++D_-)(p_-)>-1$ is satisfied. This
proves the second assertion in (a).

The proof of (a) in the case $p_+=p_-$
is similar and left to the reader.

To deduce (b) assume that $p_+=p_-=p$. As $(\alpha_+)$ is not
satisfied we have $D_+(p)+D_-(p)\ne 0$ while (\ref{alphaprim})
does not hold. By symmetry it is enough to show that the divisor
$D_\ext$ cannot be distinguished and rigid at the same time.
By Lemma \ref{lemtail}(c) either $\tO_p^-$ is not a $(-1)$-curve,
or the subtail $L_{s+2}$ is contractible. In the first case
$D_\ext$ is not rigid while in the second one it is not
distinguished. \eproof \erem

Theorem \ref{nonspec} and Remark \ref{bothnot} imply the
following.

\bcor\label{neco} (a) Under the assumptions of Theorem
\ref{nonspec} suppose additionally that $\supp\{D_+\}\cup
\supp\{D_-\}$ consists of at most one point. Then at least one of
the divisors $D_\ext, \,D_\ext^\vee$ is rigid if and only if
$(\alpha_*)$ holds. Moreover the following are equivalent:
\begin{enumerate}\item[$\bullet$] both $D_\ext, \,D_\ext^\vee$
are distinguished and rigid; \item[$\bullet$] at least one of them
is; \item[$\bullet$] $(\alpha_+)$ is fulfilled.
\end{enumerate} (b) In the case where $\supp\{D_+\}\cup
\supp\{D_-\}$ consists of two distinct points, the following are
equivalent:
\begin{enumerate}\item[$\bullet$] both $D_\ext, \,D_\ext^\vee$
are distinguished; \item[$\bullet$] at least one of them is rigid;
\item[$\bullet$] both of them are distinguished and rigid;
\item[$\bullet$] $(\beta)$ is fulfilled.
\end{enumerate}\ecor

The condition that $D_\ext, \,D_\ext^\vee$ are both distinguished
is also necessary in order that $(\beta)$ were fulfilled. Indeed,
for
$(D_+,\,D_-)=\left(-\frac{3}{2}[p_+],\,-\frac{1}{2}[p_-]\right)$
the divisor $D_\ext$ is distinguished, while $D_\ext^\vee$ is not
and both of them are non-rigid.

In the next example we exhibit two smooth  Gizatullin surfaces
completed by the same zigzag, such that one of them is a
$\C^*$-surface, whereas the second one does not admit a
$\C^*$-action, even after any logarithmic deformation keeping the
divisor at infinity fixed.

\bexa\label{dondefo} There exists a smooth Gizatullin
$\C^*$-surface, say $V_0$, with boundary zigzag
$[[0,0,-4,-2,-2]]$, see Example 4.7.3 in \cite{FKZ2}. To construct
a second Gizatullin surface, say $V$, let us consider the
following configuration $D_\ext$ in a suitable blowup $\bV\to
Q=\PP^1\times\PP^1$:

\medskip

$$D_\ext:\qquad
\cou{C_0}{0}\llin\cou{C_1}{0}\llin\cou{\quad\quad C_2}{-4}\nlin
\xbshiftup{} {\qquad\{\fF_\rho\}_{\rho= 1,2}}\vlin{30}
\cou{\hskip-8mm C_3}{-2}\nlin\cshiftup{B_2}{-1}
\llin\cou{C_4}{-2}\quad,
$$
where the map $\Phi: \bV\to Q$ is given by the linear systems
$|C_0|$ and $|C_1|$ and the feathers $\fF_1$ and $\fF_2$ consist
of two single $(-1)$-bridges. Inspecting Proposition
\ref{eboundary.9} this extended divisor $D_\ext$ does not
correspond to a Gizatullin $\C^*$-surface.

By Proposition \ref{corgen}
the divisor $D_\ext$ is stable under
generalization. However, due to Proposition \ref{convspec} it does
admit a nontrivial specialization. Namely, any of the feathers
$\fF_\rho$ can jump to $C_3$ or to $C_4$. Using Proposition
\ref{degen.3}(c) under such a specialization the dual graph of
$D_\ext$ still has at least two branching vertices and so,
 cannot correspond to a Gizatullin $\C^*$-surface,
 see Proposition \ref{eboundary.9}.

Thus indeed the surface $V=\bV\backslash D$ with
$D=C_0+\ldots+C_4$ cannot be deformed to one with a $\C^*$-action.
\eexa

\section{The reconstruction space}

Given a Gizatullin surface, any two SNC completions are related via
a birational transformation which we call a  {\em reconstruction}.
Let us denote by $\gamma$  the corresponding combinatorial
transformations of the weighted
 dual graphs of the boundary divisors.
The main result of this section (Corollary \ref{rec.8}) states that
the space of all geometric reconstructions of a pair $(X,D)$ with
a given combinatorial type $\gamma$ has a natural structure of an
affine space $\A^m$ for some $m$.

\subsection{Reconstructions of boundary zigzags}
We use in the sequel the following terminology from \cite{FKZ1}.

\bdefi\label{creconstruction} Let $\Gamma$ and $\Gamma'$ be
weighted graphs. A {\em combinatorial reconstruction} or simply
{\em reconstruction} of $\Gamma$ into $\Gamma'$ consists in a
sequence \bdi\label{bdi} \,\qquad\gamma:\quad
\Gamma=\Gamma_0&\rDotsto^{\gamma_1}& \Gamma_1 & \rDots^{\gamma_2}
&\cdots & \rDotsto^{\gamma_n}& \Gamma_n=\Gamma'\,, \edi where each
arrow $\gamma_i$ is either a blowup or a blowdown. The graph
$\Gamma'$ is called the {\em end graph} of $\gamma$. The inverse
sequence $\gamma^{-1}= (\gamma_n^{-1}, \ldots, \gamma_1^{-1})$
yields a reconstruction of $\Gamma'$ with end graph $\Gamma$.
Reconstructions can be composed: if $\gamma$ is a reconstruction
of $\Gamma$ with end graph $\Gamma'$ and $\gamma'$ is a
reconstruction of $\Gamma'$ with end graph $\Gamma''$, then the
sequence $(\gamma,\gamma')$ gives a reconstruction of $\Gamma$
into $\Gamma''$.

A reconstruction $\gamma$ is called {\em admissible} if it only
involves
\begin{enumerate}[$\bullet$]
\item blowdowns of at most linear vertices;
\item inner blowups i.e.,  blowups at edges;
\item outer blowups done at end vertices
i.e., vertices of degree $\le 1$.
\end{enumerate}
Thus an admissible reconstruction does not change the number of
branch points of the graph and their degrees.
\edefi

\bsit\label{getr}  We let $(X,D)$ and $(Y,E)$ be two pairs
consisting of smooth complete surfaces and SNC divisors on them.
Similarly as in the combinatorial setting we can speak about a
reconstruction $\tgamma$ of $(X,D)$ into $(Y,E)$ meaning a
sequence of blowups and blowdowns
$$
\bdi
\tgamma:\quad X=X_0&\rDotsto^{\tgamma_1}& X_1 &
\rDotsto^{\tgamma_2} &\cdots &
\rDotsto^{\tgamma_n}& X_n=Y\quad,
\edi
$$
performed on $D$ and on its subsequent total transforms. We say
that $\tgamma$ {\em is of  type} $\gamma$ if $\gamma$ is the
corresponding reconstruction of the dual graph $\G_D$ into $\G_E$.
Clearly the complements $X\backslash D$ and $Y\backslash E$ are
isomorphic under the birational transformation
$\tgamma:X\dasharrow Y$.

A reconstruction $\tgamma$
will be called {\em linear} if there is a domination
\bdi[small]
&& Z\\
&\ldTo &&\rdTo\\
X&& \rDashto &&Y \edi such that the total transform of $D$ is a
linear chain of rational curves.
\esit

The next fact follows immediately from Proposition 2.9 in
\cite{FKZ1}.

\bprop\label{equivariant.5} For any two standard completions
$(X,D)$ and $(Y,E)$ of a Gizatullin surface $V$ there exists an
admissible reconstruction of $(X,D)$ into $(Y,E)$. \eprop

\bprop\la{linrec} Let $\bdi\gamma:\G&\rDotsto& \G'\edi$ be an
admissible  reconstruction as in \ref{creconstruction} between two
linear chains $\G$, $\G'$, and let $D\subseteq X$ be an SNC
divisor with dual graph $\G$. Then there exists a linear
reconstruction $\bdi\tgamma:(X,D)&\rDashto & (Y,E)\edi$ of
type $\gamma$. \eprop

\bproof
Using induction on the length $n$ of $\gamma$ we may
assume that for the shorten reconstruction $\gamma':$$\,\,\,\bdi
\G=\G_0& \rDotsto &\ldots& \rDotsto& \G_{n-1}\edi$ there exists
already a linear reconstruction \bdi\tgamma': X=X_0&\rDashto
&\ldots& \rDashto X_{n-1}\edi of type $\gamma'$. Thus $X$,
$X_{n-1}$ are dominated by a blowup $Z_{n-1}$ such that the total
transform $D'$ of $D$ in $Z_{n-1}$ is linear. Since $\gamma$ is
admissible the last transform $\gamma_n$ can be either a blowdown,
an inner blowup or an outer blowup at an end vertex, see
\ref{creconstruction}.

If $\gamma_n$ is a blowdown then blowing down the corresponding
curve in $X_{n-1}$ gives a morphism $\tgamma_n: X_{n-1}\to Y$.
Obviously $\tgamma=(\tgamma',\tgamma_{n})$ is a reconstruction of
type $\gamma$ dominated by $Z:=Z_{n-1}$ and so is linear. The same
construction works in the case where $\gamma_n$ is an inner or an
outer blowup dominated by the contraction $\G_{D'}\to \G_D$.

We let $G$ denote the total transform of $D$ in $X_{n-1}$. If
$\gamma_n$ is an inner  blowup which is not dominated by the
contraction $\G_{D'}\to \G_D$ then we perform an additional blowup
$\tgamma_n:X_{n-1}\dashrightarrow Y$ at the corresponding double
point of $G$. This is dominated by the corresponding inner  blowup
$Z_{n-1}\dashrightarrow Z$. Hence $Z$ provides a linear domination
of both $X$ and $Y$, as desired.

Similarly, if $\gamma_n$ is an outer blowup at an end vertex, say,
$v_i$ of $\G_{n-1}$ which is not dominated by the contraction
$\G_{D'}\to \G_D$ then necessarily the proper transform $v_i'$ of
$v_i$ in $\G_{D'}$ is also an end vertex. In this case we perform
additionally an outer blowup $Z_{n-1}\dashrightarrow Z$ at a point
of the corresponding irreducible component $G_i'$ of $D'$ which is
not a double point of $D'$. This yields a linear domination $Z$ of
both $X$ and $Y$, as required. Now the proof is completed.
\eproof

\subsection{Symmetric reconstructions}
\bdefi\label{symme} A reconstruction of a graph $\Gamma$ is called
{\em symmetric} if it can be written in the form $(\gamma,
\gamma^{-1})$. Clearly for a symmetric reconstruction the end
graph is again $\Gamma$. \edefi

We have the following results on symmetric reconstructions.

\bprop\label{equivariant.6} \bnum[(a)] \item We let $(X,D)$ and
$(Y,E)$ be two standard completions of a normal Gizatullin surface
$V\not\cong \A^1\times \C^*$. After replacing, if necessary,
$(X,D)$ by its reversion $(X^\vee, D^\vee)$ there exists a
symmetric reconstruction of $(X,D)$ into $(Y,E)$. \item Let $X$ be
a normal surface and $D$ be a complete SNC divisor on $X$ with
dual graph $\Gamma$. Given an admissible symmetric reconstruction
$\gamma=(\tau, \tau^{-1}):\Gamma\dashrightarrow\Gamma$, there is a
reconstruction of $(X,D)$ into itself of type $\gamma$.
\enum \eprop

\bproof (a) By Proposition \ref{equivariant.5} there exists an
admissible reconstruction $\tgamma: (X,D)\dashrightarrow (Y,E)$ of
type, say, $\gamma$. Using again Proposition \ref{linrec} we can
find a linear reconstruction
$\tilde\eta:(X',D')\dashrightarrow (X,D)$ of type
$\eta:=\gamma^{-1}$, where $(X',D')$
is another standard completion of $V$.
Thus the composition
$(\tilde\eta,\tgamma):(X',D')\dashrightarrow (Y,E)$
of type $(\gamma^{-1},\gamma)$ is symmetric. We note that our
standard zigzags are different from $[[0,0,0]]$ since $V\not\cong
\A^1\times \C^*$. As follows from Proposition 3.4 in \cite{FKZ1},
any linear reconstruction
of a standard zigzag different from $[[0_{2k+1}]]$ is either the
identity or the reversion. Thus $(X',D')=(X,D)$ or
$(X',D')=(X^\vee,D^\vee)$.

(b) Clearly there is a reconstruction $\tilde\tau$ of $(X,D)$ of
type $\tau$. Then $\tgamma=(\tilde\tau,\tilde\tau^{-1})$ has the
desired properties. This completes the proof.\eproof

\subsection{Moduli space of reconstructions}
In this subsection we show that the reconstructions of a given
type form in a natural way a moduli space. \bdefi\label{rec.5} Let
$f: \cX\to S$ be a flat family of normal surfaces and
$\cD=\cD_1\cup \ldots \cup \cD_r\subseteq \cX$ be a family of SNC
divisors sitting in the smooth part of $f$. We assume that
$\cD_i\to S$ is a smooth family of curves for every $i$ and that the fiber
$\cD(s)$ forms an SNC divisor with the same dual graph $\Gamma$ in
each fiber $\cX_s$. If $\gamma$ is a reconstruction of $\Gamma$ as
in Definition \ref {creconstruction}, then a reconstruction of
$\cX/S$ of type $\gamma$ is a sequence \bdi \tgamma:\quad
\cX=\cX_0&\rDotsto^{\tgamma_1}& \cX_1 & \rDotsto^{\tgamma_2}
&\cdots & \rDotsto^{\tgamma_n}& \cX_n\quad, \edi where at each
step $\cX_{i+1}$ is either the blowup of $\cX_i$ in a section
$\iota:S\hto \cX_i$ or a blowdown of a family of $(-1)$-curves
$\cC\subseteq \cX_i$ such that fiberwise $\tgamma$ is of type
$\gamma$. \edefi

In the next result we show that the set of all reconstructions of
$X$ of type $\gamma$ has a natural structure of a smooth scheme.
It is convenient to formulate this result in a relative setup.

With the notations as in Definition \ref{rec.5},
if $S'\to S$ is a morphism of algebraic $\C$-schemes
and $\tgamma$ is a reconstruction of $\cX/S$ of
type $\gamma$ then by a base change $S'\to S$ we obtain a
reconstruction $\tgamma'$ of $\cX\times_SS'/S'$. This defines a
set valued functor $R_\gamma$ on the category of $S$-schemes that
assigns to an $S$-scheme $S'$ the set of all reconstructions of
type $\gamma$ of $\cX\times_SS'/S'$.

\bprop\label{rec.6} With $\Gamma$ and $\cX/S$ as in Definition
\ref{rec.5} the functor $R_\gamma$ is representable. The latter
means that there exists an $S$-scheme $\cR=\cR_\gamma$ of finite
type over $S$ and a universal reconstruction in $R_\gamma(\cR)$:
\bdi \tgamma_u:\quad \cX_0:=\cX\times_S\cR
&\rDotsto^{\tgamma_{u1}}& \cX_1 & \rDotsto^{\tgamma_{u2}} &\cdots
& \rDotsto^{\tgamma_{un}}& \cX_n\quad \edi such that for every
$S$-scheme $S'$ and every reconstruction $\tgamma\in R_\gamma(S')$
there is a unique $S$-morphism $g:S'\to \cR$ satisfying
$\tgamma=g^*(\tgamma_u)$. Moreover $\cR$ is smooth over $S$.
\eprop

\bproof Let us first assume   that $\gamma$ consists of a single
blowdown or an inner blowup of $\Gamma$. We claim that in these
cases $\cR:=S$ is the required moduli space. The universal family
$\tgamma_u$ is constructed as follows. If $\gamma$ is the blowdown
of the vertex corresponding to the component $\cD_\rho$ of $\cD$,
then $\cD_\rho$ is a family of $(-1)$-curves and so can be blown
down via a map $\tgamma:\cX\to \cX'$ so that $\cX'\to S$ is a flat
family, see Lemma \ref{rec.3}. It is clear that $\tgamma_u:=
\tgamma $ is in this case the universal reconstruction of type
$\gamma$.

Similarly, suppose that $\gamma$ is the blowup of the edge joining
the two vertices which correspond to $\cD_\rho$ and $\cD_\tau$. In
particular $\cD_\rho\cap \cD_\tau$ is a section of $\cX\to S$.
Blowing up this section leads to a morphism $\tau:\cX'\to \cX$,
and the composed map $\cX'\to\cX\to S$ is flat. It is easy to
check that in this case $\tgamma_u:= \tau^{-1}\in R_\gamma(S)$ is
the universal reconstruction of type $\gamma$.

We assume further that $\gamma$ is an outer blowup in a vertex of
$\Gamma$ which corresponds to $\cD_\rho$. The complement
$$
\cR:=\cD_\rho\backslash\bigcup_{\tau\ne \rho} \cD_\tau
$$
is then smooth over $S$, and the fiber product
$\cX_\cR:=\cX\times_S\cR\to \cR$ is a flat family of normal
surfaces which has a canonical section given by the diagonal
embedding $\cR\hto\cX_\cR$. The blowup $\tau: \cX'\to\cX_\cR$ of
this section provides again a universal reconstruction
$\tgamma_u:= \tau^{-1}\in R_\gamma(S)$ of type $\gamma$.

To build up the reconstruction space for an arbitrary sequence
$\gamma=(\gamma_1,\ldots, \gamma_n)$ as in Definition \ref{rec.5}
we proceed by induction on $n$. Assume that there is a universal
reconstruction space $\cR'$ for the sequence $\gamma':=(\gamma_1,
\ldots, \gamma_{n-1})$ of length $n-1$. Thus the universal
reconstruction $\tgamma'_u$ of type $\gamma'$ consists in a
sequence \bdi \tgamma_u':\quad \cX_0'=\cX\times_S
\cR'&\rDotsto^{\tgamma'_{u1}} & \cX_1' & \rDotsto^{\tgamma'_{u2}}
&\cdots & \rDotsto^{\tgamma'_{un-1}}& \cX_{n-1}' \edi as in
Definition \ref{rec.5}. Let $\cD'\subseteq \cX_{n-1}'$ be the
total transform of $\cD\times_S\cR'$ so that the dual graph of
$\cD'$ is $\Gamma_{n-1}$. Now $\gamma_n:
\Gamma_{n-1}\dashrightarrow\Gamma_n$ is a reconstruction of length
1. Hence by the first part of the proof there exists a universal
reconstruction space $\cR$ for $\cX_{n-1}'/\cR'$, where the
universal reconstruction is a birational transformation
$$
\bdi
\tgamma_{un}:\cX_{n-1}:=\cX_{n-1}'\times_{\cR'}\cR
&\rDotsto &\cX_n\,.
\edi
$$
Combining the universal  properties of
$\cR'$ and $\cR$
it follows that $\cR$ together with
$$\bdi
\tgamma_u:\quad \cX_0=\cX\times_S\cR
&\rDotsto^{\tgamma_{u1}}& \cX_1:=\cX_1'\times_{\cR'}\cR
& \rDotsto^{\tgamma_{u2}} &\cdots
& \rDotsto^{\tgamma_{un-1}}& \cX_{n-1}
&\rDotsto^{\tgamma_{un}}&\cX_n\,,
\edi
$$
where $\tgamma_{ui}:=\tgamma_{ui}' \times_{\cR'}\id_{\cR}$, forms
the required universal reconstruction of type $\gamma$.

Finally let us show that $\cR$ is smooth over $S$. Using the
iterative construction of $\cR$ it is sufficient to show this for
a reconstruction $\gamma: \Gamma\dashrightarrow \Gamma_1$ of
length 1. But the latter
 is immediate
from the first part of the proof.
\eproof

In the case where the reconstruction is admissible we get
the following important information
on the structure of $\cR$.

\bprop\label{rec.7} Let $\Gamma$, $\gamma$ and $\cX/S$ be as in
Definition \ref{rec.5}. We let $\Gamma_i$ denote the dual graph of
the total transform $\cD^{(i)}$ of $\cD$ in $\cX_i$, and we assume
that the following conditions are fulfilled: \bnum[(i)]\item
$H^1(S,\cO_S)=0$ and $\Pic(S)=0$. \item $\Gamma$ is connected, and
for every $i$ the graph $\Gamma_i$ is not reduced to a point.
\item $\gamma$ is admissible. \enum Then the reconstruction space
$\cR=\cR_\gamma$ is isomorphic to $S\times \A^m$ for some $m\in
\N$. \eprop

\bproof Let us first consider the case where the reconstruction
$\gamma:\Gamma\to \Gamma_1$ has length 1. If $\gamma$ is a
blowdown or an inner blowup we have $\cR=S$, hence the assertion
is obvious. If $\gamma$ is an outer blowup then  by our assumption
it is performed in an end vertex of $\Gamma$. The corresponding
component of $\cD$, say, $\cD_\rho$ meets exactly one other
component, say, $\cD_\tau$. The intersection $\Sigma:=\cD_\rho\cap
\cD_\tau$ is a section of the $\PP^1$-bundle $\cD_\rho\to S$. Thus
by Lemma \ref{rec.4} $\cD_\rho\to S$ is $S$-isomorphic to the
product $S\times\PP^1$ so that the section corresponds to
$S\times\{\infty\}$. Since $\cR=\cD_\rho\backslash \cD_\tau$ by
our construction, we conclude that $\cR$ is $S$-isomorphic to
$S\times\A^1$.

In the general case we proceed by induction. We consider
$\gamma'=(\gamma_1,\ldots, \gamma_{n-1})$ and the universal
reconstruction space $\cR'$ over $S$ of combinatorial type
$\gamma'$. By induction hypothesis $\cR'$ is $S$-isomorphic to
$S\times\A^{m'}$. Since $\cR=\cR_\gamma$ is the  universal
reconstruction of $\gamma_n$ with respect to\footnote{See the
proof of Proposition \ref{rec.6}.} $\cX_{n-1}\times_S\cR'/\cR'$,
from the first part of the proof we obtain that $\cR\cong \cR'$ or
$\cR\cong \cR'\times\A^1$, proving the result. \eproof

Propositions \ref{rec.6} and \ref{rec.7} lead to the following
corollary.

\bcor\label{rec.8} Let $X$ be a normal surface, and let $D$ be an
SNC divisor in $X_\reg$ with dual graph $\Gamma$. Given a
reconstruction $\gamma$ of $\Gamma$, the set $\cR_\gamma$ of all
reconstructions of $X$ of type $\gamma$ has a natural structure of
a smooth scheme. Moreover if $\gamma$ is admissible then
$\cR_\gamma\cong \A^m$ for some $m\ge 0$. \ecor

\section{Applications}

Here we prove Theorems \ref{01} and \ref{MT} on the uniqueness of
$\C^*$- and $\C_+$-actions. The proofs are based on  the results
of the previous sections and on
Theorem \ref{main} below, which states that a standard completion
of a Gizatullin surface with a distinguished and rigid extended
divisor $D_\ext$ is up to reversion (see \ref{reversion}) unique.

\subsection{The main technical result}

To formulate our result let us first fix the notations.

\bsit\la{unisit1} Let $V$ be a non-toric Gizatullin surface and let
$(\bV, D)$ and $(\bV', D')$ be standard completions of $V$.
We also consider
the minimal resolutions of singularities $V'$, $(\tV,D)$,
$(\tV',D')$ of $V$, $(\bV,D)$  and $(\bV', D')$, respectively.
As in \ref{unisit} we let
$$
\Phi=\Phi_0\times\Phi_1:\tV\to \PP^1\times\PP^1
\quad\mbox{and}\quad
\Phi'=\Phi'_0\times\Phi'_1:\tV'\to \PP^1\times\PP^1
$$
denote the standard morphism
and $D_\ext$,  $D'_\ext$ the  extended divisors.

Reversing the zigzag $D'=[[0,0,w'_2,\ldots, w'_n]]$ by a sequence
of inner elementary transformations provides the standard
completion $(\bV^{\prime\vee}, D^{\prime\vee})$, see
\ref{reversion}. \esit

\bthm\label{main}
Assume that the extended divisor $D_\ext$ of $(\tV, D)$ is
distinguished and rigid.
After replacing $(\bV', D')$ by $(\bV^{\prime\vee}, D^{\prime\vee})$
if necessary there is an isomorphism $f:\bV\to \bV'$
with $f(D)=D'$.
\ethm

Note that this isomorphism is {\em not} the identity
on the affine part $V$, in general.

\bproof
Replacing $(\bV', D')$ by $(\bV^{\prime\vee}, D^{\prime\vee})$
if necessary, by Proposition \ref{equivariant.6}(a)
there is a reconstruction
$\tgamma'$ from $(\tV,D)$ to $(\tV',D')$ of
type, say $\gamma$, which is admissible and symmetric.
Thus $\tgamma'$ can be
considered as a point in the reconstruction space
$\cR=\cR_\gamma\cong\A^m$, see Corollary \ref{rec.8}.
By Proposition \ref{equivariant.6}(b) there is also a
reconstruction $\tgamma$ of
$(X,D)$ of type $\gamma$ into itself. Let
$$
\bdi \tgamma_u:\quad \cX_0=\tV\times\cR &\rDotsto^{\tgamma_{u1}}&
\cX_1 & \rDotsto^{\tgamma_{u2}} &\cdots &
\rDotsto^{\tgamma_{un-1}}& \cX_{n-1}
&\rDotsto^{\tgamma_{un}}&\cX_n\, \edi
$$
be the universal reconstruction of combinatorial
type $\gamma$ and consider the family
$\tilde\cV:=\cX_n$ together with the total transform
$\cD$ of $D\times\cR$ in $\tilde\cV$.
Thus $(\tilde\cV,\cD)$ is a family of completions
of $V'$ over the reconstruction space $\cR$  as considered in
Proposition \ref{propdist}. Moreover, by construction the completions
$(\tV, D)$ and $(\tV', D')$ are the fibers over the points $\tgamma$,
$\tgamma'\in \cR$, respectively.

Let now $\cD_\ext$ be the family of extended divisors of $(\tilde
\cV, \cD)$. Its fibre over $\tgamma$ is $D_\ext$ and so is rigid.
Hence the family of extended divisors $\cD_\ext$ has the same dual
graph over each point of $\cR$.  By Proposition \ref{propdist} the
family $(\tilde \cV, \cD)$ is trivial and so there is an
isomorphism $ (\tV, D)\times \cR\cong (\tilde \cV, \cD)$.
Restricting it to the fibre over $\tgamma'$ gives an isomorphism
$\tilde f:(\tV, D)\to (\tV',D')$ that induces an isomorphism
$f:\bV\to \bV'$ with the desired property. \eproof

In particular, in the situation of Theorem \ref{main} it follows
that the extended divisors $D_\ext$, $D'_\ext$, considered as
schemes via their reduced structures, are isomorphic at least
after reversion, if necessary. It is important to note that  this holds even without the
assumption that $D_\ext$ is distinguished:.

\bprop\label{mainprop} With the notations as in \ref{unisit1},
assume that the extended divisor $D_\ext$ of $(\tV, D)$ is rigid.
After replacing $(\bV', D')$ by $(\bV^{\prime\vee},
D^{\prime\vee})$, if necessary, the corresponding extended
divisors are isomorphic as reduced curves under an isomorphism
$\tilde f:D_\ext\to D'_\ext$ with $\tilde f(D)=D'$ preserving the
weights. \eprop

\bproof As in the proof above the family of extended divisors
$\cD_\ext$ has the same dual graph over each point of $\cR\cong
\A^m$.
Since the fibers of $\cD_\ext$ are trees of rational curves
with at least 2 components, the result is immediate from Corollary
\ref{rec.51}.
\eproof

\subsection{Uniqueness of $\C^*$-actions}
In Theorem \ref{main*} below we deduce part (1) of Theorem
\ref{MT}.

\bthm\label{main*}
Let $V=\Spec \C[t][D_+,D_-]$ be a non-toric normal Gizatullin surface
satisfying one of the following two conditions.
\begin{enumerate}
\item[($\alpha_*$)] $\supp \{D_+\}\cup \supp \{D_-\}$ is empty or
consists of one point, say $p$, where
$$
D_+(p) + D_-(p) \le -1
\quad\mbox{or}\quad \{D_+(p)\} \ne 0\ne \{ D_-(p)\}.
$$
\item[($\beta$)] $\supp \{D_+\}=\{p_+\}$ and $\supp
\{D_-\}=\{p_-\}$ for two distinct points $p_+,p_-$, where
$$
D_+(p_+) + D_-(p_+) \le -1 \quad\mbox{and}\quad D_+(p_-) + D_-(p_-)
\le -1\,.
$$
\end{enumerate}
Then the $\C^*$-action on $V$ is unique, up to conjugation in the
group $\Aut (V)$ and up to inversion $\lambda\mapsto \lambda^{-1}$
in $\C^*$. Moreover the given $\C^*$-action is conjugate to its
inverse if and only if there is an automorphism $\psi:\A^1\to\A^1$
such that
$$
\psi^*(D_-)-D_+ \mbox{ is integral and } \psi^*(D_++D_-)=D_++D_-\;.
$$
\ethm

\bproof Let $\Lambda, \Lambda':\C^*\times V\to V$ be two
$\C^*$-actions on $V$, where $\Lambda$ is the given one. We
consider the corresponding equivariant standard completions
$(\bV,D)$ and $(\bV',D')$ of $V$. After reversing the first one,
if necessary, its extended divisor $D_\ext$ is rigid according to
Theorem \ref{nonspec}. Applying Proposition \ref{mainprop}, after
reversing $(\bV',D')$, if necessary, the extended divisors
$D_\ext$ and $D'_\ext$ are isomorphic. Since by
Proposition 5.12 in \cite{FKZ2} and its proof a non-toric
Gizatullin $\C^*$-surface is uniquely determined by its extended
divisor,  the first part follows. The second one is a
consequence of Lemma \ref{symmetric0}.
\eproof

Applying Theorem \ref{main*} to smooth Gizatullin $\C^*$-surfaces,
we obtain the following

\bcor\label{smooth*} If $V=\Spec\C[t][D_+,D_-]$ is a non-toric
smooth Gizatullin $\C^*$-surface, then its $\C^*$-action is
uniquely determined up to conjugation and inversion unless its
standard zigzag is \be\label{00w} [[0,0, (-2)_{s-2}, w_{s},
(-2)_{n-s}]]\,, \quad\mbox{where} \quad w_s\le -2,\,\,n\ge 4
\mbox{ and } 2\le s\le n\,. \ee \ecor

\bproof
Suppose that $\supp \{D_+\}\subseteq \{p_+\}$ and $\supp
\{D_-\}\subseteq \{p_-\}$. If $p_+=p_-=:p$ and $\{D_+(p)\}\ne
0\ne \{D_-(p)\}$ then by Theorem \ref{main*} the $\C^*$-action is
unique up to conjugation and inversion. Otherwise either $p_+\ne
p_-$ or one of the fractional parts $\{D_+\}$, $\{D_-\}$ vanishes.
Anyhow the smoothness of $V$ implies the desired form (\ref{00w})
of the dual graph of $D$, see Remark \ref{vosst}(2).
\eproof

\subsection{Uniqueness of $\A^1$-fibrations}
\bsit \label{5.6} In this subsection  we consider a normal
Gizatullin surface $V$ with a fixed standard completion $(\bV,
D)$, where $D=[[0,0,w_2,\ldots,w_n]]$  is a zigzag with
irreducible components $C_0,\ldots, C_n$. As usual the linear
system $|C_0|$ defines an $\A^1$-fibration $\Phi_0:V\to \A^1$.
Reversion as in \ref{reversion}  provides the standard completion
$(\bV^\vee, D^\vee)$ so that $D^\vee$ has irreducible components
$C^\vee_0,\ldots, C^\vee_n$ with self-intersections
$[[0,0,w_n,\ldots,w_2]]$. The linear system $|C^\vee_0|$ defines a
second $\A^1$-fibration $\Phi^\vee_0:V\to \A^1$, which we call the
{\em reverse fibration}. We say that two $\A^1$-fibrations $\varphi,
\varphi':V\to\A^1$ are {\em conjugate} if
$\varphi'=\beta\circ\varphi\circ\alpha$ for some automorphisms
$\alpha$ of $V$ and $\beta$ of $\A^1$.
\esit

In Theorem \ref{main+} below we give a partial answer to the
following problem.

\begin{conj}
{\it Suppose that $V$ is not a Danilov-Gizatullin surface. Is then
every $\A^1$-fibration $\varphi:V\to \A^1$ conjugate to one of the
standard $\A^1$-fibrations $\Phi_0, \Phi_0^\vee$?}
\end{conj}

The latter is actually equivalent to the uniqueness problem
for $\C_+$-actions on $V$ in the sense of (3) and (4) below.
Let us recall some standard facts concerning $\C_+$-actions.

\bsit\label{Ren}
1. (\cite{Ren}) If $\C_+$ acts on
an affine algebraic $\C$-scheme $V=\Spec A$
then the associated derivation $\partial$ on $A$ is locally
nilpotent, i.e.\ for every $f\in A$ we can find $n\in \N$ such
that $\partial^n(f)=0$. Conversely, given a locally nilpotent
$\C$-linear derivation $\partial:A\to A$ the map
$\varphi:\C_+\times A\to A$ with $\varphi(t,f):=e^{t\partial}f $
defines an action of $\C_+$ on $V$.

2.  (See e.g., \cite{ML1, Zai}) Assume that $A$ as in (1) is a
domain and let $\p\in {\Der}_{\C} A$ be a locally nilpotent
derivation of $A$. Then the subalgebra $\ker \p=A^{\C_+}\subseteq
A$ is algebraically and factorially closed, or inert\footnote{The
latter means that $ab\in\ker \p\Rightarrow a,\,b\in\ker \p$.}, in
$A$, and the field extension $\Frac (\ker \p) \subseteq \Frac A$
has transcendence degree $1$. Moreover for any $u\in\Frac A$ with
$u\p(A)\subseteq A$, the derivation $u\p\in \Der_{\C} A$ is
locally nilpotent if and only if $u\in\Frac(\ker\p)$.

If $A$ as in (1) is normal then the ring of invariants $A^{\C_+}$
is normal too. If  $\dim A\le 3$ then by a classical result of
Zariski \cite{Za} $A^{\C_+}$ is finitely generated and
$C=\Spec\,A^{\C_+}$  is the algebraic quotient $V//\C_+$. Thus the
orbit map $V\to C$ provides an $\A^1$-fibration.

3. Conversely if a normal affine surface $V$
admits an $\A^1$-fibration $V\to C$ over a smooth affine curve
$C$, then there exists a non-trivial regular $\C_+$-action on $V$
along this fibration. It is unique up to multiplication of an
infinitesimal generator $\partial$ with an element $u\in
\Frac(\ker\p)$ as in (2).

4. As mentioned in the introduction, every normal affine surface
$V$ which is not a Gizatullin surface admits
at most one $\A^1$-fibration over $\A^1$, see \cite{BML}.
\esit

We restrict in the sequel to $\A^1$-fibrations on Gizatullin
surfaces. Let us provide several examples of such
fibrations.

\bexa \label{fib*} 1. Let $V=\C[t][D_+,D_-]$ be a Gizatullin
$\C^*$-surface. Taking in \ref{5.6} an equivariant standard
completion the $\A^1$-fibrations $\Phi_0$, $\Phi_0^\vee$ on $V$
are equivariant with respect to suitable $\C^*$-actions on $\A^1$.
By Proposition 3.25 in \cite{FlZa2}, they are given by two
homogeneous elements \be\label{+/-fibration} v_+:V\to \A^1
\quad\mbox{and}\quad v_-:V\to \A^1 \ee of positive and negative
degree, respectively. Moreover, by {\it loc. cit.} any other
$\A^1$-fibration $\varphi:V\to\A^1$ compatible with the
$\C^*$-action on $V$ is equal to $v_+$ or $v_-$.

2. The toric surface $V_{d,e}=\A^2\quot \Z_d$ (see \ref{toricsit})
admits many hyperbolic $\C^*$-actions. Indeed, for any coprime
integers $a,b$ the action $t.(x,y):=(t^ax,t^by)$, $t\in \C^*$, on
$\A^2$ descends to $V$, and in the case where $ab<0$ it is
hyperbolic. Up to a twist, the $\A^1$-fibrations $v_\pm:V\to\A^1$
are induced by the projections $(x,y)\mapsto x$, $(x,y)\mapsto y$,
respectively.

3. Let now $V=V_{k+1}$ be a Danilov-Gizatullin surface, see \cite
{FKZ2}, section 5.3. According to {\it loc.cit.}, Corollary
5.16(b) $V$ carries at least $ \lfloor \frac{k+1}{2}\rfloor$
pairwise non-conjugate $\A^1$-fibrations $V_{k+1}\to\A^1$. \eexa

The following theorem is the main result of this subsection.

\bthm\label{main+} Let $V$ be a Gizatullin surface with a
distinguished and rigid extended divisor\footnote{This is
fulfilled for instance if the assumptions of Theorem
\ref{maincrit} hold.} $D_\ext$. Then every $\A^1$-fibration
$\varphi: V\to\A^1$ is conjugate to one of $ \Phi_0$,
$\Phi_0^\vee$. \ethm

Before starting the proof, let us make the following observation.

\bsit\label{sest} Consider  a semistandard completion\footnote{See
\ref{sst}.} $(\bV', D')$  of a Gizatullin surface $V$, where
$D'=C'_0+\ldots+C'_n$ and $(C_0')^2=0$. Then the linear system
$|C'_0|$ defines a morphism $\Phi_0':\bV'\to \PP^1$ which
restricts to an $\A^1$-fibration $V\to\A^1$.

Conversely, we claim that {\it any $\A^1$-fibration
$\varphi:V\to \A^1$ is induced by the standard $\A^1$-fibration of
a suitable standard completion $(\bV, D)$ of $V$.} Indeed,
given an $\A^1$-fibration $\varphi:V\to \A^1$, there exists an
effective $\C_+$-action on $V$ along this fibration, see
\ref{Ren}(3). By virtue of Lemma \ref{equivariant}(c) one can find an equivariant semistandard completion $(\bV',
D')$ of $V$ such that $\varphi$ extends to a morphism
$\varphi':\bV'\to \PP^1$. Performing a sequence of elementary
transformations with centers at the fiber $C_0'$ of $\varphi'$,
one can reach a standard completion, say, $(\bV, D)$ of $V$, where
this time $D=C_0+\ldots+C_n$ with $C_0^2=C_1^2=0$. The morphism
$\Phi_0:\bV\to\PP^1$ defined by the linear system $|C_0|$
restricts again to $\varphi:V\to \A^1$.
\esit

\bproof[Proof of Theorem \ref{main+}.] We let as in \ref{5.6}
$(\bV, D)$ denote the standard completion of $V$ with standard
$\A^1$-fibration $\Phi_0$, and we let $(\bV',D')$ denote another
such standard pair with standard morphism  as in \ref{sest}
inducing the given fibration $\varphi:V\to \A^1$.

Since by our assumption the extended divisor $D_\ext$ is
distinguished and rigid, Theorem \ref{main} applies. By this
theorem, $(\bV, D)$ is isomorphic to one of $(\bV', D')$,
$(\bV^{\prime\vee}, D^{\prime\vee})$ or, equivalently, $(\bV',
D')$ is isomorphic to one of the pairs $(\bV', D')$,
$(\bV^{\prime\vee}, D^{\prime\vee})$.  In particular $\varphi$ is
conjugate to $\Phi_0$ or $\Phi_0^\vee$ under this
isomorphism. \eproof

The following lemma shows that the extended divisor is uniquely
determined by $\varphi$.

\blem\la{rem+} Let $(\bV,D)$ and $(\bV',D')$ be two standard
completions  of the same Gizatullin surface $V$. If the associated
$\A^1$-fibrations $\Phi_0, \Phi_0':V\to \A^1$ are conjugate  then there is an isomorphism
$f:D_\ext\to D'_\ext$ of the corresponding extended divisors
(regarded as reduced curves) with
$f(D)=D'$, which preserves the weights.
\elem

\bproof
We may assume that the automorphism of $V$ which
conjugates $\Phi_0$ and $\Phi_0'$ extends to a birational map
$\bdi \tilde f:\tV& \rDashto &\tV'\edi$ of the minimal resolutions
of $\bV,\,\bV'$ with $\Phi_0'\circ\tilde f=\Phi_0$. If
$D=C_0+\ldots+ C_n$ and $D'=C_0'+\ldots+ C_n'$ then clearly
$\tilde f$ is regular at the points of $C_1\backslash (C_0\cup
C_2)$. Performing elementary transformations on $\tV$ with centers
at $C_0$, if necessary, we may suppose that $\tilde f$ is
biregular along $C_0$, so that ${\tilde f}^{-1}$ is also regular
along $(C_0'\cup C_1')\backslash C_2$. Contracting the
divisors\footnote{Both of them have negatively definite
intersection forms.} $C_2+\ldots+ C_n$ and $C_2'+\ldots+ C_n'$ on
the surfaces $\tV$ and $\tV'$ to singular points $p$, $p'$,
respectively, yields two normal surfaces $W$ and $W'$. Moreover
$\tilde f$ induces a birational map $\bar f: W\to W'$ which is an
isomorphism outside $p,p'$. By the Riemann extension theorem $\bar
f$ is actually an isomorphism. Then also $\tilde f$, obtained from
$\bar f$ via minimal resolution of singularities, is an
isomorphism. Hence $\tilde f$ induces an isomorphism of the
boundaries and the extended divisors of the two completions. Since
$(C_1')^2=0$, also $C_1^2=0$ and so the standard zigzag $D$
remains the same under the above elementary transformations. Now
the lemma follows. \eproof

Let us apply these results to a $\C^*$-surface $V=\Spec
\C[t][D_+,D_-]$. In this case we choose in \ref{5.6} the
equivariant standard completion $(\bV,D)$ so that $\Phi_0$ and
$\Phi_0^\vee$ are equivariant. The next result yields part (2) of
Theorem \ref{MT}.

\bcor\label{cor+*} We let $V=\Spec \C[t][D_+,D_-]$ be a Gizatullin
$\C^*$-surface. If one of the conditions ($\alpha_+$), ($\beta$)
of \ref{greek} is fulfilled,  then the following hold. \bnum[(1)]
\item Every $\A^1$-fibration $V\to\A^1$ is conjugate to one of
$\Phi_0$ or $\Phi_0^\vee$. \item Assume furthermore that $V$ is
non-toric. The $\A^1$-fibrations $\Phi_0$, $\Phi^\vee_0$ are then
conjugate if and only if $\{D_+(p_+)\}=\{D_-(p_-)\}$
and the divisor $D_++D_-$ is stable under an automorphism of
$\A^1$ interchanging $p_+$ and $p_-$. In the latter case up to
conjugation there is only one $\A^1$-fibration $V\to\A^1$. \enum
\ecor

\bproof By Theorem \ref{nonspec} under our assumptions the
extended divisor $D_\ext$ is distinguished and rigid. So (1)
follows directly from Theorem \ref{main+}. To deduce (2), assume
first that  $\{D_+(p_+)\}=\{D_-(p_-)\}$ and
$D_++D_-=\psi^*(D_++D_-)$ for an appropriate automorphism
$\psi\in\Aut (\A^1)$ interchanging $p_+$ and $p_-$. By Lemma
\ref{symmetric0} the $\C^*$-surfaces $\Spec A_0[D_+,D_-]$ and
$\Spec A_0[D_-,D_+]$ with $A_0=\C[t]$ are isomorphic. This
isomorphism interchanges the fibrations $v_+$ and $\psi\circ v_-$
 as in Example \ref{fib*}(1). Hence $\Phi_0$, $\Phi^\vee_0$ are
conjugate.

Suppose now that $\Phi_0$, $\Phi^\vee_0$ are conjugate. By Lemma \ref{rem+}  there is an isomorphism of extended
divisors $f:D_\ext\to D'_\ext$ as reduced curves  with $f(D)=D'$
preserving the weights. According to Proposition 5.12 in
\cite{FKZ2} and its proof the $\C^*$-surfaces $\Spec A_0[D_+,D_-]$
and $\Spec A_0[D_-,D_+]$ are equivariantly isomorphic. Now the
assertion follows from Lemma \ref{symmetric0}. \eproof

As a particular case we obtain the following result, which was
proved in the smooth case by Daigle \cite{Dai} and Makar-Limanov
\cite{ML2}.

\bcor Let $V$ be a normal surface in $\A^3_\C$ with equation
$xy=P(t)$, where $P(t)\ne 0$ is a polynomial. Then every
$\A^1$-fibration on $V$ is conjugate to $x:V\to \A^1$. \ecor

\bproof According to Example 4.10 in \cite{FlZa2}, $V$ admits a
DPD presentation $V=\Spec \C[t][D_+,D_-]$ with integral divisors
$D_+=0$ and $D_-=-\div (P)$. Thus condition $(\alpha_+)$ is
fulfilled and so the result follows from Corollary
\ref{cor+*}(1,2) in virtue of Remark \ref{srt}. \eproof

Let us finally examine $\A^1$-fibrations of affine toric surfaces.

\bprop\la{toric+} The toric surface $V_{d,e}\cong \A^2\quot\Z_d$
(see \ref{toricsit}) admits at most 2 conjugacy classes of
$\A^1$-fibrations over $\A^1$. Moreover, there is only one such
conjugacy class if and only if $e^2\equiv 1\mod d$. \eprop

\bproof The DPD presentation of $V_{d,e}$ considered in the proof of
Lemma \ref{toric} satisfies $(\alpha_+)$.
Applying Corollary \ref{cor+*} gives the first part. To prove the
second assertion, we assume first that
$e^2\equiv 1\mod d$. Using the notations of Example \ref{fib*}(1),
(2) the affine fibrations $\Phi_0$, $\Phi_0^\vee$ are induced by
the projections $(x,y)\mapsto x$ and $(x,y)\mapsto y$. Because of
our assumption the map $h: (x,y)\mapsto (y,x)$ satisfies
$h(\zeta.(x,y))=\zeta^e.(y,x)$. Hence $h$ induces an automorphism
$\bar h$ on the quotient $V_{d,e}$ that interchanges these
projections and thus also $\Phi_0$ and $\Phi_0^\vee$.

Conversely assume that the $\A^1$-fibrations $\Phi_0$,
$\Phi_0^\vee$ are conjugate in $\Aut(V)$. According to Lemma
\ref{rem+} the standard zigzag $D$ of $V$ is symmetric. Due to
Lemma \ref{toric} $D$ and the reversed zigzag $D^\vee$ are given
by
$$
D:\qquad\quad\co{0}\lin\co{0}\llin \boxo{
\frac{d-e}{d}}\quad,\qquad D^\vee:\qquad\quad\co{0}\lin\co{0}\llin
\boxo{ \frac{d-e'}{d}}\quad,
$$ where $0\le e,e'<d$ and $ee'\equiv 1\mod d$, cf. \ref{not1}.
Hence $D$ and $D^\vee$ are equal if and only if $e=e'$ or,
equivalently, $e^2\equiv 1\mod d$. \eproof

\end{document}